\documentclass[twoside,11pt]{article}
\title{} \author{} \date{}
\usepackage{latexsym,amssymb,times}%,showkeys}
\input amssym.def
\newtheorem{te}{Theorem}[section]
\newtheorem{fac}[te]{Fact}
\newtheorem{lem}[te]{Lemma}
\newtheorem{cla}[te]{Claim}
\newtheorem{rem}[te]{Remark}

\def\dok{\noindent{\bf Proof. }}
\def\kdok{\hfill $\Box$ \par \vspace*{2mm} }
\def\dom{\mathop{\mathrm{dom}}\nolimits}
\def\ran{\mathop{\mathrm{ran}}\nolimits}
\def\Emb{\mathop{\mathrm{Emb}}\nolimits}
\def\Age{\mathop{\mathrm{Age}}\nolimits}
\def\dense{\mathop{\mathrm{Dense}}\nolimits}
\def\supp{\mathop{\mathrm{supp}}\nolimits}
\begin{document}
\thispagestyle{plain}
\begin{center}
           {\large \bf
           MAXIMAL CHAINS OF ISOMORPHIC SUBORDERS OF  \\[1mm]
                    COUNTABLE ULTRAHOMOGENEOUS PARTIAL ORDERS}

\end{center}
\begin{center}
{\bf Milo\v s S.\ Kurili\'c\footnote{Department of Mathematics and Informatics, University of Novi Sad,
     Trg Dositeja Obradovi\'ca 4, 21000 Novi Sad, Serbia. e-mail: milos@dmi.uns.ac.rs}
     and
     Bori\v sa  Kuzeljevi\'c\footnote{Mathematical Institute of the Serbian Academy of Sciences and Arts,
     Kneza Mihaila 36, 11001 Belgrade, Serbia. e-mail: borisa@mi.sanu.ac.rs}
}
\end{center}
\begin{abstract}
\noindent
We investigate the poset  $\langle \mathbb P (\mathbb X )\cup \{\emptyset \}, \subset  \rangle$, 
where $\mathbb P (\mathbb X )$ is the set of isomorphic suborders of a countable
ultrahomogeneous partial order $\mathbb X$. For $\mathbb X$
different from (resp.\ equal to) a countable antichain the order
types of maximal chains in $\langle \mathbb P (\mathbb X )\cup \{\emptyset \}, \subset \rangle$ 
are characterized as the order types of compact (resp.\ compact and nowhere dense) sets of reals
having the minimum non-isolated.

\noindent
{\sl 2000 Mathematics Subject Classification}:
06A06, % Partial order,general
06A05, % Total order
03C15, % Denumerable structures
03C50. % Models with special properties.
\\
{\sl Keywords}: ultrahomogeneous partial order, isomorphic substructure, maximal chain, compact set.
\end{abstract}
\section{Introduction}\label{S1}
The general concept - to explore the relationship between the
properties of a relational structure $\mathbb X$ and the
properties of the poset $\mathbb P (\mathbb X )$ of its isomorphic
substructures - can be developed in several ways.
For example, regarding the forcing theoretic aspect, the
poset of copies of each countable non-scattered linear order is
forcing equivalent to the two-step iteration of the Sacks forcing
and a $\sigma$-closed forcing \cite{KurTod}, while the posets of
copies of countable scattered linear orders have $\sigma$-closed
forcing equivalents (separative quotients) \cite{Kscatt}.

Regarding the order-theoretic aspect, one of extensively
investigated order invariants of a poset is the class of order
types of its maximal chains \cite{Day,Kopp,Kura,Monk} and, for the
poset of isomorphic suborders of the rational line, $\langle
\mathbb Q , <_\mathbb Q\rangle$, this class is characterized in
\cite{K1}. The main result of the present paper is the following
generalization of that result.
\begin{te} \rm \label{T5032}
If $\mathbb X $ is a countable ultrahomogeneous partial order
different from a countable antichain, then for each linear order
$L$ the following conditions are equivalent:

(a) $L$ is isomorphic to a maximal chain in the poset $\langle
\mathbb P (\mathbb X ) \cup \{ \emptyset \} , \subset \rangle $;

(b) $L$ is an ${\mathbb R}$-embeddable complete linear order with $0_L$ non-isolated;

(c) $L$ is isomorphic to a compact set $K \subset \mathbb R$ such
that $0_K \in K^\prime$.

\noindent If $\mathbb X $ is a countable antichain, then the
corresponding characterization is obtained if we replace
``complete" by ``Boolean" in (b) and ``compact" by ``compact and
nowhere dense" in (c).
\end{te}
So, for example, there are maximal chains of copies of the random poset isomorphic to $(0,1]$, to the Cantor set without 0, and to $\alpha ^*$, for each countable limit ordinal $\alpha$.
Although it is not a usual practice, we start with a proof in the introduction. The equivalence of (b) and (c) is a known fact
(see, for example, Theorem 6 of \cite{K1}) and the implication (a) $\Rightarrow$ (b) follows from the general result on ultrahomogeneous
structures given in Theorem \ref{T5037} of the present paper. Thus, only the implication (b) $\Rightarrow$ (a) remains to be proved. Naturally, we will use the following,
well known classification of countable ultrahomogeneous partial orders - the Schmerl list \cite{Sch}:
\begin{te}[Schmerl]\rm\label{T4027}
A countable strict partial order is ultrahomogeneous iff it is isomorphic to one of the following partial orders:

 $\mathbb A _\omega$, a countable antichain (that is, the empty relation on $\omega$);

 $\mathbb B _n = n \times \mathbb Q$, for $1\leq n\leq \omega$, where $\langle i_1, q_1\rangle < \langle i_2 ,q_2\rangle \Leftrightarrow i_1=i_2 \;\land \; q_1 <_\mathbb Q q_2$;

 $\mathbb C _n=n\times \mathbb Q$, for $1\leq n\leq \omega$, where  $\langle i_1, q_1\rangle < \langle i_2 ,q_2\rangle \Leftrightarrow  q_1 <_\mathbb Q q_2$;

 $\mathbb D$, the unique countable homogeneous universal poset (the random poset).
\end{te}
For the antichain $\mathbb A _\omega$ the implication (b)
$\Rightarrow$ (a) follows from Theorem \ref{T5023} and the fact
that $\mathbb P(\mathbb A _\omega)=[\omega ]^\omega$ is a positive
family. The most difficult part of the proof of (b) $\Rightarrow$
(a) - for the random poset $\mathbb D$ - is given in Section
\ref{S5}. In Sections \ref{S6} and \ref{S7}, using the
constructions from \cite{K1}, we prove  (b) $\Rightarrow$ (a) for
the posets $\mathbb B _n$ and $\mathbb C _n$.

The rest of this section contains two facts which will be used in the sequel.
We remind the reader that a linear order $\langle L,<\rangle $ is called
{\it Boolean} iff it is {\it complete} (has 0,1 and has no gaps) and {\it has dense jumps}, which means that for each $x,y \in L$
satisfying $x<y$ there are $a,b \in L$ such that $x\leq a< b \leq y$ and $(a,b)_L=\emptyset$.
\begin{fac}  \rm \label{T5022}
Each countable complete linear order is Boolean.
\end{fac}
We recall that a family ${\mathcal P} \subset P(\omega)$ is called a {\it positive family} iff:

(P1) $\emptyset \notin {\mathcal P}$;

(P2) ${\mathcal P} \ni A \subset B \subset \omega \Rightarrow  B \in {\mathcal P}$;

(P3) $A \in {\mathcal P} \wedge |F| <\omega \Rightarrow  A \backslash F \in {\mathcal P}$;

(P4) $\exists A \in{\mathcal P} \;\; |\omega \backslash A|=\omega$.
\begin{te}  \rm \label{T5023}
(\cite{K})
If ${\mathcal P}\subset P(\omega)$ is a positive family, then for each linear order $L$ the
following conditions are equivalent:

(a) $L$ is isomorphic to a maximal chain in the poset $\langle
{\mathcal P} \cup\{\emptyset\},\subset \rangle$;

(b) $L$ is an $\mathbb R$-embeddable Boolean linear order with
$0_L$ non-isolated;

(c) $L$ is isomorphic to a compact nowhere dense set $K\subset
\mathbb R$ such that $0_K\in K'$.

\noindent
In addition, (b) implies that there is a maximal chain $\mathcal L$ in
$\langle {\mathcal P} \cup \{\emptyset\}, \subset \rangle$ satisfying $\bigcap (\mathcal L \setminus \{ \emptyset \}) = \emptyset$ and isomorphic to $L$.
\end{te}
\section{Copies of countable ultrahomogeneous structures}\label{S3}
Let $L=\{ R_i :i\in I \}$ be a relational language, where
ar$(R_i)=n_i$, $i\in I$. An $L$-structure $\mathbb X = \langle X,
\{ \rho _i :i\in I \} \rangle$ is called {\it countable} iff
$|X|=\omega$. If $A\subset X$, then $\langle A, \{ (\rho _i)_A
:i\in I \} \rangle$ (shortly denoted by $\langle A, \{ \rho _i
:i\in I \} \rangle$, whenever this abuse of notation does not
produce a confusion) is a {\it substructure} of $\mathbb X$, where
$(\rho_i)_A = \rho _i \cap A^{n_i}$, $i\in I$. If $\mathbb Y
=\langle Y, \{ \sigma _i : i\in I\} \rangle$ is an $L$-structure
too, a mapping $f:X \rightarrow Y$ is an {\it embedding} (we write
$\mathbb X \hookrightarrow _f \mathbb Y$) iff it is an injection
and
$$\forall i\in I \;\;\forall \langle x_1, \dots   x_{n _i}\rangle \in X^{n_i} \;\;
(\langle x_1, \dots   ,x_{n _i}\rangle \in \rho _i \Leftrightarrow
\langle f(x_1), \dots  ,f(x_{n _i})\rangle \in \sigma _i).
$$
If $\mathbb X$ embeds in $\mathbb Y$ we write $\mathbb X
\hookrightarrow \mathbb Y$. Let $\Emb (\mathbb X , \mathbb Y )  =
\{ f: \mathbb X \hookrightarrow _f \mathbb Y \} $ and $\Emb
(\mathbb X ) = \{ f: \mathbb X \hookrightarrow _f \mathbb X \}$.
If, in addition, $f$ is a surjection, it is an {\it isomorphism}
(we write $\mathbb X \cong _f \mathbb Y$) and the structures
$\mathbb X$ and $\mathbb Y$ are {\it isomorphic}, in notation
$\mathbb X \cong \mathbb Y$.

A {\it finite isomorphism} of $\mathbb X$ is each isomorphism
between finite substructures of $\mathbb X$. A structure $\mathbb
X$ is {\it ultrahomogeneous} iff each finite isomorphism on
$\mathbb X$ can be extended to an automorphism of $\mathbb X$. The
{\it age} of $\mathbb X$, $\Age \mathbb X$, is the class of all
finite $L$-structures embeddable in $\mathbb X$. We will use the
following well known facts from the Fra\"{\i}ss\'{e} theory.
\begin{te}[Fra\"{\i}ss\'{e}] \rm \label{T5036}
Let $L$ be an at most countable relational language. Then

(a) A countable $L$-structure $\mathbb X$ is ultrahomogeneous iff
for each finite isomorphism $\varphi$ of $\mathbb X$  and each
$x\in X\setminus \dom \varphi$ there is a finite isomorphism
$\psi$ of $\mathbb X$ extending $\varphi$ to $x$ (see \cite{Fra}
p.\ 389 or \cite{Hodg} p.\ 326).

(b) If $\mathbb X$ and $\mathbb Y$ are countable ultrahomogeneous
$L$-structures and $\Age \mathbb X =\Age \mathbb Y$, then $\mathbb
X \cong \mathbb Y$ (see \cite{Fra} p.\ 333 or \cite{Hodg} p.\
326).
\end{te}
Concerning the order types of maximal chains in the posets of the
form $\langle \mathbb P (\mathbb X ), \subset \rangle$, where
$\mathbb X=\langle X, \{ \rho _i :i\in I \} \rangle$ is a
relational structure and $\mathbb P (\mathbb X )$ the set of the
domains of its isomorphic substructures, that is
$$
\mathbb P (\mathbb X )  =  \{ A\subset X : \langle A, \{ (\rho
_i)_A :i\in I \} \rangle \cong \mathbb X\}
          =  \{ f[X] : f \in \Emb (\mathbb X )\}
$$
we have the following general statement.
\begin{te} \rm \label{T5037}
Let $\mathbb X$ be a countable ultrahomogeneous structure of an at
most countable relational language and $\mathbb P (\mathbb X )\neq
\{ X\}$. If $\mathcal L$ is a maximal chain in the poset  $\langle
\mathbb P (\mathbb X )\cup \{ \emptyset \} , \subset \rangle$,
then

(a) $\mathcal L$ is an $\mathbb R$-embeddable complete linear
order with $0_\mathcal L (=\emptyset)$ non-isolated;

(b) If there is  a positive family $\mathcal P \subset \mathbb P
(\mathbb X )$, then for each countable linear order $L$ satisfying
(a), there is a maximal chain in $\langle \mathbb P (\mathbb X
)\cup \{ \emptyset \} , \subset \rangle$  isomorphic to $L$.
\end{te}
\dok
(a) First we prove that
\begin{equation}\label{EQ5021}\textstyle
\bigcup \mathcal A \in \mathbb P (\mathbb X ), \mbox{ for each
chain } \mathcal A \mbox{ in the poset } \langle \mathbb P
(\mathbb X ) , \subset \rangle .
\end{equation}
Let $\varphi$ be a finite isomorphism of $\bigcup \mathcal A $ and
$x\in \bigcup \mathcal A$. Since $\mathcal A$ is a chain there is
$A\in \mathcal A$ such that $\dom \varphi \cup \ran \varphi \cup
\{ x \}\subset A$. Since $A\cong \mathbb X$, by Theorem
\ref{T5036}(a) there is $y\in A$ such that $\psi = \varphi \cup \{
\langle x,y \rangle\}$ is an isomorphism so $\psi$ is
a finite isomorphism of $\bigcup \mathcal A$. Thus, by Theorem
\ref{T5036}(a), the structure $\bigcup \mathcal A$ is
ultrahomogeneous. Since $\mathbb X \cong A \subset \bigcup
\mathcal A \subset X$ we have $\Age \mathbb X =\Age A \subset \Age
\bigcup \mathcal A \subset \Age \mathbb X$, which, by Theorem
\ref{T5036}(b), implies $\bigcup \mathcal A \cong \mathbb X$, that
is $\bigcup \mathcal A \in \mathbb P (\mathbb X )$.

Let $X=\{x_n:n \in \omega\}$ be an enumeration. Since ${\mathcal
L} \subset [X]^{\omega} \cup \{ \emptyset \}$, the function
$f:{\mathcal L}\to {\mathbb R}$ defined by $f(A)=\sum _{n \in
\omega } 2^{-n}\cdot \chi_A (x_n)$ (where $\chi_A :X \to \{0,1\}$
is the characteristic function of the set $A \subset X$) is an
embedding of $\langle{\mathcal L},\subset\rangle$ into $\langle
{\mathbb R}, <_\mathbb R \rangle$.

Clearly, $\min {\mathcal  L}=\emptyset$ and $\max {\mathcal L}=X$.
Let $\langle{\mathcal A}, {\mathcal B}\rangle$ be a cut in
${\mathcal L}$. If ${\mathcal A}=\{ \emptyset \}$ then $\max
{\mathcal A}= \emptyset$. If ${\mathcal A} \neq \{ \emptyset \}$,
by (\ref{EQ5021}) we have $\bigcup{\mathcal  A}\in \mathbb P
(\mathbb X )$ and, since $A\subset \bigcup{\mathcal  A}\subset B$,
for each $A\in {\mathcal  A}$ and $B\in {\mathcal  B}$, the
maximality of ${\mathcal  L}$ implies $\bigcup{\mathcal  A} \in
{\mathcal  L}$. So, if $\bigcup{\mathcal  A}\in {\mathcal  A}$
then $\max {\mathcal A}=\bigcup{\mathcal  A}$. Otherwise
$\bigcup{\mathcal A}\in {\mathcal B}$ and $\min{\mathcal  B}=
\bigcup{\mathcal A}$. Thus $\langle{\mathcal L},\subset\rangle$ is
complete.

Suppose that $A$ is the successor of $\emptyset$ in ${\mathcal
L}$. Since $\mathbb P (\mathbb X )\neq \{ X\}$ there is $B\in
\mathbb P (\mathbb X)\setminus \{ X \}$ and, if $f:\mathbb X
\hookrightarrow A$, then $f[B]\in \mathbb P (\mathbb X )$,
$f[B]\varsubsetneq A$ and, hence, $\mathcal L \cup \{ f[B] \}$ is
a chain in $\mathbb P (\mathbb X )$. A contradiction to the
maximality of ${\mathcal  L}$.

(b) By Fact \ref{T5022}, $L$ is a Boolean order and, by Theorem
\ref{T5023}, in the poset $\langle \mathcal P \cup \{ \emptyset\},
\subset \rangle$ there is a maximal chain $\mathcal L$ isomorphic
to $L$ and such that $\bigcap (\mathcal L \setminus \{ \emptyset
\} ) = \emptyset$. Now, ${\mathcal L}$ is a chain in $\langle
\mathbb P (\mathbb X ) \cup \{ \emptyset \} , \subset \rangle$ and
we check its maximality. Suppose that ${\mathcal L} \cup \{ A \} $
is a chain, where $A\in \mathbb P (\mathbb X ) \setminus
{\mathcal L}$. Then $A\varsubsetneq S$ or $S\varsubsetneq A$, for
each $S\in {\mathcal L}\setminus \{ \emptyset \}$ and, since
$\bigcap ({\mathcal L}\setminus \{ \emptyset \} )=\emptyset$,
there is $S\in {\mathcal L}\setminus \{ \emptyset \}$ such that
$S\subset A$, which implies $A\in \mathcal P$. But ${\mathcal
L}\setminus \{ \emptyset \}$ is a maximal chain in $\mathcal P$. A
contradiction. \hfill $\Box$
\begin{rem}  \rm \label{R5000}
Concerning the assumption $\mathbb P (\mathbb X )\neq \{ X\}$ we
note that there are countable ultrahomogeneous structures
satisfying $\mathbb P (\mathbb X )= \{ X\}$ (see \cite{Fra}, p.\
399).

For $1<n<\omega$ the set $\mathbb P(\mathbb C_n)$ does not contain
a positive family, since (P3) is not satisfied. Namely, if
$A\in \mathbb P(\mathbb C_n)$ and $x\in A$, then $A\setminus
\{x\}$ is not a copy of $\mathbb C_n$ (one class of incompatible
elements is of size $n-1$).

For some $\omega$-saturated, $\omega$-homogeneous-universal
relational structures the implication (b) $\Rightarrow$ (a) of
Theorem \ref{T5032} is not true. Let $L$ be the language with one
binary relational symbol $\rho$ and ${\mathcal T}$ the $L$-theory
of empty relations ($\forall x,y \;\neg\; x\,\rho\; y$). Then
$\mathbb X=\langle \omega,\emptyset \rangle$ is the
$\omega$-saturated model of ${\mathcal T}$. But $\mathbb P(\mathbb
X)=[\omega]^{\omega}$ is a positive family and, by Theorem
\ref{T5023}, maximal chains in $\mathbb P(\mathbb X)\cup
\{\emptyset \}$ are Boolean. Thus, for example, $\mathbb P(\mathbb
X)\cup \{\emptyset \}$ does not contain a maximal chain isomorphic
to $[0,1]_\mathbb R$.
\end{rem}
\section{Copies of the countable random poset}\label{S4}
%\begin{df}\rm\label{D5000}
Let $\mathbb P= \langle P, < \rangle$ be a partial order. By
$C(\mathbb P )$ we denote the set of all triples $\langle L,G,U
\rangle$ of pairwise disjoint finite subsets of $P$ such that:

(C1) $\forall l\in L \; \forall g\in G \;\; l<g$,

(C2) $\forall u\in U \; \forall l\in L \;\; \neg u<l $ and

(C3) $\forall u\in U \; \forall g\in G \;\; \neg g<u $.

\noindent For $\langle L,G,U \rangle \in C(\mathbb P )$, let $P
_{\langle L,G,U \rangle}$ be the set of all $p\in P\setminus (L
\cup G\cup U)$ satisfying:

(S1) $\forall l\in L \; p>l $,

(S2) $\forall g\in G \; p<g $ and

(S3) $\forall u\in U \; p \| u $ (where $p\| q$ denotes that $p\neq q \land \neg p< q \land \neg q< p$).
%\end{df}
\begin{fac}\rm\label{T5003}
Let $\mathbb P= \langle P, < \rangle$ be a partial order and
$\emptyset \neq A\subset P$. Then

(a) $C(A, < )= \{ \langle L,G,U \rangle \in C(\mathbb P ): L,G,U
\subset A\}$;

(b) $A_{\langle L,G,U \rangle}= P_{\langle L,G,U \rangle}\cap A$,
for each $\langle L,G,U \rangle \in C(A, < )$.

(c) $\langle \emptyset , \emptyset , \emptyset \rangle \in
C(\mathbb P )$ and $P_{\langle \emptyset , \emptyset , \emptyset
\rangle}=P$.
\end{fac}
\dok For pairwise disjoint sets  $L,G,U\subset A$ we have:
$L\times G \subset \;<$ iff $L\times G \subset \;< _A$ and
$((U\times L)\cup (G\times U))\;\cap < \; =\emptyset $ iff $((U\times L)\cup (G\times U))\;\cap < _A \; =\emptyset $.
\hfill $\Box$
\begin{fac}\rm\label{T5040}
A countable strict partial order $\mathbb D = \langle D, <
\rangle$ is a countable random poset iff $D_{\langle L, G ,U
\rangle}\neq \emptyset$, for each $\langle L, G ,U \rangle \in
C(\mathbb D )$ (see \cite{Cam}).
\end{fac}
\begin{lem}\rm\label{T5006}
Let $\mathbb D = \langle D, < \rangle$ be a countable random
poset. Then

(a) $D_{\langle L, G ,U \rangle}\in \mathbb P (\mathbb D )$ and,
hence, $|D_{\langle L, G ,U \rangle}|=\omega$, for each $\langle
L, G ,U \rangle \in C(\mathbb D )$;

(b) $D\setminus F \in \mathbb P (\mathbb D )$, for each finite
$F\subset D$;

(c) If $D = A\;\dot{\cup}\; B$, then either $A$ or $B$ contains an
element of $\mathbb P (\mathbb D )$;

(d) If $\mathcal L\subset \mathbb P (\mathbb D)$ is a chain, then
$\bigcup\mathcal L\in \mathbb P (\mathbb D)$;

(e) If $C\subset D$ and $A\not\subset C$ for each $A\in \mathbb P
(\mathbb D)$, then $D\setminus C\in \mathbb P (\mathbb D )$.
\end{lem}
\dok (a) Let $\langle L, G ,U \rangle \in C(\mathbb D )$. Then
$L$, $G$ and $U$ are disjoint subsets of $D$,
\begin{equation}\label{EQ5018}
\forall l\in L \;\; \forall g\in G \;\;\forall u\in U \;\;(u\not< l < g \not< u ),
\end{equation}
and $D_{\langle L, G ,U \rangle}\cap (L\cup G\cup U)=\emptyset$.
Let $\langle L_1, G_1 ,U_1 \rangle \in C(D_{\langle L, G ,U
\rangle})$. Then $L_1$, $G_1$ and $U_1$ are disjoint subsets of
$D_{\langle L, G ,U \rangle}$ and, by Fact \ref{T5003}, $\langle
L_1, G_1 ,U_1 \rangle \in C(\mathbb D )$ which implies
\begin{equation}\label{EQ5019}
\forall l_1 \in L_1 \;\; \forall g_1 \in G_1 \;\;\forall u_1\in U_1 \;\;(u_1 \not< l_1 < g_1 \not< u_1 ).
\end{equation}
Since $L_1\cup G_1 \cup U_1 \subset D_{\langle L, G ,U \rangle}$,
by (S1)-(S3) we have
\begin{equation}\label{EQ5020}
\forall x \in L_1\cup G_1 \cup U_1 \;\; \forall l\in L \;\; \forall g\in G \;\;\forall u\in U \;\;(l < x < g \land x\not< u \land u\not< x ).
\end{equation}
First we show that $\langle L\cup L_1,G\cup G_1, U\cup U_1
\rangle\in C(\mathbb D)$. (C1) Let $l'\in L\cup L_1$ and $g'\in
G\cup G_1$. Then $l' < g'$ follows from: (\ref{EQ5018}), if $l'\in
L$ and $g'\in G$; (\ref{EQ5019}), if $l'\in L_1$ and $g'\in G_1$;
(\ref{EQ5020}), if $l'\in L$ and $g'=x\in G_1$ or $l'=x\in L_1$
and $g'\in G$. (C2) Let $l'\in L\cup L_1$ and $u'\in U\cup U_1$.
Then $u' \not< l'$ follows from: (\ref{EQ5018}), if $l'\in L$ and
$u'\in U$; (\ref{EQ5019}), if $l'\in L_1$ and $u'\in U_1$;
(\ref{EQ5020}), if $l'\in L$ and $u'=x\in U_1$ (since $l'<u'$) or
$l'=x\in L_1$ and $u'\in U$. In the same way we prove (C3).

So there is $x\in D_{\langle L\cup L_1,G\cup G_1, U\cup U_1
\rangle}$, which implies $x\in D_{\langle L, G ,U \rangle}\cap
D_{\langle L_1, G_1 ,U_1 \rangle}$ $  =(D_{\langle L, G ,U
\rangle})_{\langle L_1, G_1 ,U_1 \rangle}$ (Fact \ref{T5003}).
Thus $D_{\langle L, G ,U \rangle}$ is a random poset and, hence a
copy of $\mathbb D$.

(b) Let $\langle L,G,U \rangle\in C(D\setminus F)$. By Fact
\ref{T5003} we have $\langle L,G,U \rangle\in C(\mathbb D )$ and,
by (a), $\emptyset\neq (D\setminus F) \cap D_{\langle L, G ,U
\rangle}= (D\setminus F)_{\langle L, G ,U \rangle}$. Thus
$D\setminus F$ is a copy of $\mathbb D$.

(c) Suppose that $P(A)\cap \mathbb P (\mathbb D)=\emptyset$. Then
$A\notin \mathbb P (\mathbb D)$ and, hence, there is $\langle
L,G,U \rangle\in C(A)$ such that $A_{\langle
L,G,U\rangle}=D_{\langle L, G ,U \rangle}\cap A =\emptyset$. By
Fact \ref{T5003} we have $\langle L,G,U \rangle\in C(\mathbb D )$
and, by (a), $\mathbb P (\mathbb D )\ni D_{\langle L, G ,U
\rangle}\subset B$.

(d) See (\ref{EQ5021}) in the proof of Theorem  \ref{T5037}.

(e) Let $\langle L,G,U \rangle\in C(D\setminus C)$. Then, by Fact
\ref{T5003}, $\langle L,G,U \rangle\in C(\mathbb D )$ and, by (a),
$D_{\langle L, G ,U \rangle}\in \mathbb P (\mathbb D )$. By the
assumption we have $D_{\langle L, G ,U \rangle}\cap (D\setminus
C)\neq\emptyset$ and, by Fact \ref{T5003}, $(D\setminus
C)_{\langle L, G ,U \rangle}\neq \emptyset$ and $D\setminus C$ is
a random poset. \hfill $\Box$
\begin{lem}\rm\label{T5018}
Let $\mathbb D = \langle D, < \rangle$ be a countable random
poset, $C\in [D]^\omega$ and $A\not\subset C$ for each $A\in
\mathbb P (\mathbb D )$ (for example, $C$ can be an infinite
antichain). Then

(a) $\mathcal P = \{ B\subset D : D\setminus C \subset ^* B\}
\subset \mathbb P (\mathbb D ) $;

(b) $\mathcal P$ is a positive family on $D$.
\end{lem}
\dok (a) Suppose that $A\subset D\setminus B$, for some $A\in
\mathbb P (\mathbb D )$. Since $D\setminus C \subset ^* B$ we have
$D\setminus B \subset ^* C$ and, hence, $A\subset ^* C$, that is
$|A\setminus C|<\omega$. By Lemma \ref{T5006}(b), $A\cap
C=A\setminus (A\setminus C)\in \mathbb P (\mathbb D )$, which is
not true. So $D\setminus B$ does not contain copies of $\mathbb D$
and, by Lemma \ref{T5006}(e), $B\in \mathbb P (\mathbb D )$.

(b) Conditions (P1) and (P2) are evident. If $D\setminus C \subset
^* B$ and $|F|<\omega$, then, clearly,  $D\setminus C \subset ^*
B\setminus F$ and (P3) is true. Since the set $D\setminus C$ is
co-infinite (P4) is true as well. \hfill $\Box$
\begin{lem}  \rm \label{T5028}
Let $A\subset B\subset \omega$ and let $L$ be a complete linear
ordering, such that $|B\setminus A|=|L|-1$. Then there is a chain
${\mathcal L}$ in $[A,B]_{P (B )}$ satisfying $A,B\in {\mathcal
L}\cong L$ and such that $\bigcup {\mathcal  A},\bigcap {\mathcal
B}\in {\mathcal  L}$ and $|\bigcap {\mathcal  B}\setminus \bigcup
{\mathcal  A}|\leq 1$, for each cut $\langle {\mathcal
A},{\mathcal B}\rangle$  in ${\mathcal L}$.
\end{lem}
\dok
If $|B\setminus A|$ is a finite set, say $B=A \cup \{ a_1 , \dots a_n \}$,
then $|L|=n+1$ and $\mathcal L =\{ A, A\cup \{ a_1 \}, A\cup \{ a_1 , a_2  \} , \dots , B\}$ is a chain
with the desired properties.

If $|B\setminus A|=\omega$, then $L$ is a countable and, hence, $\mathbb R$-embeddable complete linear order.
It is known that an infinite linear order is isomorphic to a maximal chain in $P(\omega)$ iff
it is ${\mathbb R}$-embeddable and Boolean (see, for example, \cite{K}).
By Fact \ref{T5023} $L$ is a Boolean order and, thus, there is a maximal chain $\mathcal L _1$ in $P(B\setminus A)$
isomorphic to $L$. Let $\mathcal L =\{ A\cup C : C\in \mathcal L _1\}$. Since $\emptyset , B\setminus A \in \mathcal L_1$
we have $A,B \in \mathcal L$ and the function $f: \mathcal L _1 \rightarrow \mathcal L$, defined by $f(C)=A\cup C$,
witnesses that $\langle \mathcal L _1 , \varsubsetneq \rangle \cong \langle \mathcal L , \varsubsetneq \rangle$ so
$\mathcal L$ is isomorphic to $L$. For each cut $\langle \mathcal A , \mathcal B\rangle$ in $\mathcal L_1$
we have $\bigcup \mathcal A \subset \bigcap \mathcal B$ and, by the maximality of $\mathcal L_1$,
$\bigcup \mathcal A , \bigcap \mathcal B \in \mathcal L _1$  and
$| \bigcap \mathcal B \setminus \bigcup \mathcal A | \leq 1$. Clearly, the same is true for
each cut in $\mathcal L$.
\hfill $\Box$
\section{Maximal chains of copies of the random poset} \label{S5}
\begin{te}  \rm \label{T5041}
For each $\mathbb R$-embeddable complete linear order $L$ with
$0_L$ non-isolated there is a maximal chain in $\langle \mathbb
P(\mathbb D)\cup \{\emptyset\},\subset \rangle$ isomorphic to $L$.
\end{te}
\dok
By Lemma \ref{T5018} and Theorem \ref{T5037} it remains to prove the statement for uncountable $L$'s.
So let $L$ be an uncountable linear order with the given properties.
\begin{cla}  \rm \label{T5024}
$L\cong \sum_{x \in [-\infty , \infty ]} L_x$, where

(L1) $L_x$, $x\in [-\infty , \infty ]$, are at most countable complete linear orders,

(L2) The set $M=\{x \in [-\infty , \infty ]:|L_x| >1 \}$ is at most countable,

(L3) $|L_{-\infty }|=1$ or $0_{L_{-\infty }}$ is non-isolated.
\end{cla}
\dok
$L=\sum _{i\in I}L_i$, where $L_i$ are the equivalence classes corresponding to the condensation relation $\sim$ on $L$ given by:
$x\sim y \Leftrightarrow |[\min \{ x,y \} ,\max \{ x,y \} ]|\leq \omega$ (see \cite{Rosen}).
Since $L$ is complete and ${\mathbb R}$-embeddable $I$ is too and,
since the cofinalities and coinitialities of $L_i$'s are countable, $I$ is a dense linear order; so $I\cong [0,1]\cong [-\infty , \infty ]$. Hence $L_i$'s are complete and, since $\min L_i \sim \max L_i$, countable.
If $|L_i|>1$,  $L_i$ has a jump (Fact \ref{T5022}) so, $L\hookrightarrow {\mathbb R}$ gives $|M|\leq \omega$.
\kdok

\noindent  {\bf Case I:} $-\infty \not\in M\ni \infty$. First we
take the rational line $\langle \mathbb Q , <_\mathbb Q \rangle$
and construct a set $\lhd \subset \mathbb Q^2$ such that $\langle
\mathbb Q,\lhd \rangle$ is a random poset with additional,
convenient properties. Let $\mathbb P$ be the set of pairs
$p=\langle P_p , \vartriangleleft _p\rangle$ satisfying

(i) $P_p \in [\mathbb Q ]^{<\omega }$,

(ii) $\vartriangleleft _p \subset P_p \times P_p$ is a strict partial order on $P_p$,

(iii) $< _\mathbb Q$ extends $\vartriangleleft _p$, that is
$\forall q_1,q_2 \in P_p \;\; (q_1 \vartriangleleft _p q_2
\Rightarrow q_1 <_\mathbb Q q_2)$,

\noindent and let the relation $\leq$ on $\mathbb P$ be defined
by:
\begin{equation} \label{EQ5001}
p\leq q \Leftrightarrow P_p \supset P_q \; \land \vartriangleleft _p \cap (P_q \times P_q)=\vartriangleleft _q.
\end{equation}
\begin{cla}\rm\label{T5008}
$\langle \mathbb P , \leq \rangle$ is a partial order.
\end{cla}
\dok
The reflexivity of $\leq $ is obvious. If $p\leq q\leq p$, then $P_p=P_q$ and, hence,
$\lhd_p =\lhd_p\cap (P_p\times P_p)=\lhd_p\cap (P_q\times P_q)=\lhd_q $ so $p=q$ and $\leq$ is antisymmetric.

If $p\leq q \leq r$, then $P_p\supset P_q \supset P_r$ and, consequently,
$\lhd_p\cap (P_r\times P_r)=\lhd_p\cap (P_q\times P_q )\cap (P_r\times P_r)= \lhd_q\cap (P_r\times P_r)=\lhd_r$. Thus $p\leq r$.
\hfill $\Box$
\begin{cla}\rm\label{T5009}
The sets $\mathcal D _q =\{ p\in \mathbb P : q\in P_p\}$, $q\in
\mathbb Q$, are dense in $\mathbb P$.
\end{cla}
\dok If $p\in \mathbb P\setminus \mathcal D_q$, that is $q\notin
P_p$, then $\lhd _p$ is an irreflexive and transitive relation on
the set $P_p$ and on the set $P_p\cup\{q\}$ as well. Also $\lhd _p
\subset <_ \mathbb Q$ thus $p_1=\langle P_p\cup\{q\},\lhd_p
\rangle \in \mathbb P$. Thus $p_1\in \mathcal D_q$ and, clearly,
$p_1\leq p$.
\kdok
\noindent Let $\mathbb Q = J \cup \bigcup
_{y\in M}J_y$ be a partition of $\mathbb Q$ into $|M|+1$ dense
subsets of $\mathbb Q$. For $\langle L,G,U \rangle \in ([\mathbb Q
]^{<\omega })^3 \setminus \{ \langle \emptyset, \emptyset,
\emptyset \rangle \}$, let $m_{\langle L,G,U \rangle}=\max
_{\langle \mathbb Q , <_\mathbb Q \rangle} (L\cup G\cup U)$.
\begin{cla}\rm\label{T5010}
For each $\langle L,G,U \rangle \in ([\mathbb Q ]^{<\omega
})^3\setminus \{ \langle \emptyset, \emptyset, \emptyset \rangle
\}$ and each $m\in \mathbb N$ the set $\mathcal D _{\langle L,G,U
\rangle , m}$ is dense in $\mathbb P$, where
\begin{eqnarray*}\textstyle
\mathcal D _{\langle L,G,U \rangle , m}&=&\Big\{ p\in \mathbb P : L\cup G\cup U\subset P_p\ \wedge\ \Big(\langle L,G,U \rangle \not\in C(p)\\
                       & &\lor \;\;\textstyle (G\neq \emptyset \land p_{\langle L,G,U \rangle}\cap J\neq \emptyset )\\
                       & &\lor \;\;\textstyle (G=\emptyset \land p_{\langle L,G,U \rangle}\cap (m_{\langle L,G,U \rangle} , m_{\langle L,G,U \rangle} +\frac{1}{m})\cap J\neq \emptyset)\Big)\Big\}.
\end{eqnarray*}
\end{cla}
\dok Let $p'\in \mathbb P\setminus\mathcal D_{\langle L,G,U
\rangle, m}$. By Claim \ref{T5009} there is $p\in \mathbb P$ such
that $p\leq p'$ and $L\cup G\cup P\subset P_p$. If $\langle L,G,U
\rangle \not\in C(p)$ then $p\in \mathcal D _{\langle L,G,U
\rangle , m}$ and we are done. If
\begin{equation}\label{EQ5007}
\langle L,G,U \rangle \in C(p) ,
\end{equation}
then we continue the proof distinguishing the following two cases.

{\it Case 1:} $G\neq \emptyset$. Let us define $\max
_{\langle\mathbb Q , <_\mathbb Q \rangle}\emptyset =-\infty$. By
(\ref{EQ5007}) and (C1) for $p$, if $L\neq \emptyset$, then
$\max_{\langle\mathbb Q , <_\mathbb Q \rangle} L \lhd _p
\min_{\langle\mathbb Q , <_\mathbb Q \rangle} G$ and, by (iii),
$\max_{\langle\mathbb Q , <_\mathbb Q \rangle} L <_\mathbb Q
\min_{\langle\mathbb Q , <_\mathbb Q \rangle} G$. Now, since $J$
is a dense set in $\langle \mathbb Q , <_\mathbb Q \rangle$ we
choose
\begin{equation}\label{EQ5008}\textstyle
q\in (\max_{\langle\mathbb Q , <_\mathbb Q \rangle} L ,
\min_{\langle\mathbb Q , <_\mathbb Q \rangle} G) \cap J\setminus
P_p
\end{equation}
and define $p_1=\langle P_p\cup\{q\},\lhd_{p_1} \rangle $ where
\begin{equation}\label{EQ5009}\textstyle
\lhd_{p_1}=\lhd_{p}\cup
\{\langle x,q \rangle : \exists l\in L\ \; x\trianglelefteq _p l\}\cup
\{\langle q,y \rangle : \exists g\in G\ \; g\trianglelefteq_p y\}.
\end{equation}
First we prove that $p_1\in \mathbb P$. Clearly, $p_1$ satisfies condition (i).

(ii) Since $\lhd _p$ is an irreflexive relation and, by (\ref{EQ5008}), $q\not\in P_p$, by (\ref{EQ5009}) the relation $\lhd_{p_1}$ is irreflexive as well.

Suppose that $\lhd_{p_1}$ is not asymmetric. Then, since  $\lhd_p$
is asymmetric, there is $t\in P_p$ such that $\langle t,q \rangle
, \langle q,t \rangle\in \lhd_{p_1}$ and by (\ref{EQ5009}),
$g\trianglelefteq _p t \trianglelefteq _p l$, for some $l\in L$
and $g\in G$ which, by the transitivity of $\trianglelefteq _p$
implies $g\trianglelefteq _p l$. But, by (\ref{EQ5007}) and (C1)
we have $l\lhd _p g$. A contradiction.

Let $\langle a,b \rangle, \langle b,c \rangle\in \lhd_{p_1}$. Then, since the relation $\lhd_{p_1}$ is irreflexive and asymmetric, we have
$a\neq b \neq c \neq a$. If $q\not\in\{ a,b,c\}$, then $\langle a,c \rangle\in\lhd_{p_1}$ by the transitivity of $\lhd_p$. Otherwise we have three possibilities:

$a=q$. Then $ \langle b,c \rangle\in \lhd_p$ and there is $g\in G$ such that $g\trianglelefteq _p b$. Hence $g\lhd _p c$ which, by
(\ref{EQ5009}), implies $\langle q,c \rangle\in \lhd_{p_1}$, that is $\langle a,c \rangle\in \lhd_{p_1}$.

$b=q$.
Then there are $l\in L $ and $g\in G$  such that $a\trianglelefteq _p l$ and  $g\trianglelefteq _p c$. By (C1) we have $l\lhd _p g$ and,
by the transitivity of $\lhd _p$, $a\lhd _p c$ and, hence, $\langle a,c \rangle\in \lhd_{p_1}$.

$c=q$. Then  $ \langle a,b \rangle\in \lhd_p$ and there is $l\in L$ such that $b\trianglelefteq _p l$.
Hence $a\lhd _p l$ which, by (\ref{EQ5009}), implies $\langle a,q \rangle\in \lhd_{p_1}$, that is $\langle a,c \rangle\in \lhd_{p_1}$.

(iii) Since $p\in \mathbb P$, we have $\lhd _p \subset <_\mathbb
Q$. If $\langle x,q \rangle \in \lhd_{p_1}$ and $l\in L$, where $x
\trianglelefteq _p l$, then, since  $\lhd _p$ satisfies (iii), we
have $x\leq _\mathbb Q l$. By (\ref{EQ5008}) we have $l<_\mathbb Q
q$ and, thus, $x<_\mathbb Q q$. In a similar way we show that
$\langle q,y \rangle \in \lhd_{p_1}$ implies $q<_\mathbb Q y$.

Thus $p_1\in \mathbb P $, $P_{p_1}\supset P_p \supset L\cup G\cup
U$ and, by  (\ref{EQ5009}), $ \lhd_{p_1} \cap(P_p \times P_p)=\lhd
_p$, which implies that $p_1 \leq p \;(\leq p')$. So $p$ is a
suborder of $p_1$ and, by (\ref{EQ5007}) and Fact \ref{T5003},
$\langle L,G,U \rangle \in C(p_1)$. Since $G\neq \emptyset$ and
$q\in J$, for a proof that $p_1\in \mathcal D _{\langle L,G,U
\rangle , m}$ it remains to be shown that $q\in (p_1)_{\langle
L,G,U\rangle}$. By (\ref{EQ5009}) $l\lhd_{p_1} q\lhd_{p_1} g$, for
each $l\in L$ and $g\in G$, so (S1) and (S2) are true. For $u\in
U$, $\langle u, q\rangle \in \lhd_{p_1}$ would give $l\in L$
satisfying $u\trianglelefteq _p l$ and, since $U\cap L=\emptyset$,
$u\lhd _p l$, which is impossible by (\ref{EQ5007}) and (C2).
Similarly, $\langle q, u\rangle \in \lhd_{p_1}$ is not possible
and, thus, $q
\parallel _{p_1}u$ and (S3) is satisfied.

{\it Case 2:} $G=\emptyset$. Again, since $J$ is a dense set in
the linear order $\langle \mathbb Q , <_\mathbb Q \rangle$ we
choose
\begin{equation}\label{EQ5010}\textstyle
q\in (m_{\langle L,G,U \rangle},m_{\langle L,G,U \rangle}+\frac{1}{m})\cap J \setminus P_p
\end{equation}
and define $p_1=\langle P_p\cup\{q\},\lhd_{p_1} \rangle $, where
\begin{equation}\label{EQ5011}\textstyle
\lhd_{p_1}=\lhd_{p}\cup
\{\langle x,q \rangle : \exists l\in L\ \; x\trianglelefteq _p l\}.
\end{equation}
First we prove that $p_1\in \mathbb P$. Clearly, $p_1$ satisfies condition (i).

(ii) By (\ref{EQ5010}) we have $q\not\in P_p$ so, by (\ref{EQ5011}) the relation $\lhd_{p_1}$ is irreflexive.

Let $\langle a,b \rangle, \langle b,c \rangle\in \lhd_{p_1}$. If
$q\not\in\{ a,b,c\}$, then $\langle a,c \rangle\in\lhd_{p_1}$ by
(\ref{EQ5011}) and the transitivity of $\lhd_p$. Otherwise, by
(\ref{EQ5011}) again, $a,b\neq q$ and, thus, $c=q$. Hence there is
$l\in L$ such that $b \trianglelefteq _p l$. Since $a,b\neq q$, by
(\ref{EQ5011}) we have $a\lhd _p b$ and, hence $a\lhd _p l$, which
implies $\langle a,q \rangle \in \lhd_{p_1}$, that is  $\langle
a,c \rangle \in \lhd_{p_1}$.

(iii)  Since $p\in \mathbb P$, we have $\lhd _p \subset <_\mathbb
Q$. If $\langle x,q \rangle \in \lhd_{p_1}$ and $l\in L$, where $x
\trianglelefteq _p l$, then, since  $\lhd _p$ satisfies (iii), we
have $x\leq _\mathbb Q l$. By (\ref{EQ5010}) we have $l\leq
_\mathbb Q m_{\langle L,G,U \rangle}<_\mathbb Q q$ and, thus,
$x<_\mathbb Q q$.

Thus $p_1\in \mathbb P $. As in Case 1 we show that $ L\cup G\cup
U \subset P_{p_1}$, $p_1 \leq p \;(\leq p')$ and $\langle L,G,U
\rangle \in C(p_1)$. By (\ref{EQ5010}) and since $G= \emptyset$,
for a proof that $p_1\in \mathcal D _{\langle L,G,U \rangle , m}$
it remains to be shown that $q\in (p_1)_{\langle L,G,U\rangle}$.
(S2) is trivial and, by (\ref{EQ5011}), for $l\in L$ we have
$\langle l,q \rangle \in \lhd_{p_1}$ thus (S1) holds as well.
Suppose that $\neg \; q
\parallel _{p_1} u$, for some $u\in U$. Then, by (\ref{EQ5010})
and (\ref{EQ5011}), $\langle u,q \rangle \in  \lhd_{p_1}$ and,
hence, there is $l\in L$ satisfying $u\lhd _p l$, which is
impossible by (\ref{EQ5007}) and (C2) for $p$. So (S3) is true.
\kdok
\noindent By the Rasiowa Sikorski theorem there is a filter
$\mathcal G $ in $\langle \mathbb P , \leq \rangle$ intersecting
the sets $\mathcal D _q$, $q\in \mathbb Q$, and $\mathcal D
_{\langle L,G,U \rangle , m}$, $\langle L,G,U \rangle \in
([\mathbb Q ]^{<\omega })^3$, $m\in \mathbb N$.
\begin{cla}\rm\label{T5034}
(a) $\bigcup _{p\in \mathcal G }P_p=\mathbb Q$;

(b) $\vartriangleleft = \bigcup _{p\in \mathcal G}\vartriangleleft
_p$ is a strict partial order on $\mathbb Q$;

(c) $ \lhd \cap (P_p \times P_p)=\lhd _p $, for each $p\in
\mathcal G $;

(d) $< _\mathbb Q$ extends $\vartriangleleft $, that is $\forall
q_1,q_2 \in \mathbb Q \; (q_1 \vartriangleleft  q_2 \Rightarrow
q_1 <_\mathbb Q q_2)$.
\end{cla}
\dok (a) For $q\in \mathbb Q$ let $p_0\in \mathcal G \cap \mathcal
D_q$. Then $q\in P_{p_0} \subset \bigcup_{p\in \mathcal G}P_p$.

(b) The relation $\lhd$ is irreflexive since all the relations $\lhd_p$ are irreflexive.

Let $\langle a,b \rangle,\langle b,c \rangle\in\lhd$, $\langle a,b
\rangle\in\lhd_{p_1}$ and $\langle b,c \rangle \in \lhd_{p_2}$,
where $p_1,p_2\in \mathcal G$. Since $\mathcal G $ is a filter
there is $p\in \mathcal G$ such that $p\leq p_1 ,p_2 $, which by
(\ref{EQ5001}) implies $\lhd _{p_1}, \lhd _{p_2}\subset \lhd _p$.
Thus $\langle a,b \rangle , \langle b,c \rangle \in \lhd_p$ and,
by the transitivity of $\lhd_p $, $\langle a,c \rangle \in \lhd_p
\subset \lhd$.

(c) The inclusion ``$\supset$" follows from (ii) and the
definition of $\lhd$. If $\langle a,b \rangle \in \lhd \cap (P_p
\times P_p)$, then there is $p_1\in \mathcal G$ such that $\langle
a,b \rangle \in \lhd _{p_1}$ and, since $\mathcal G $ is a filter,
there is $p_2\in \mathcal G$ such that $p_2\leq p ,p_1 $. By
(\ref{EQ5001}) we have $\lhd _{p_1}\subset \lhd _{p_2}$, which
implies $\langle a,b \rangle \in \lhd _{p_2}$ and, by
(\ref{EQ5001}) again, $\langle a,b \rangle \in \lhd _{p_2}\cap
(P_p \times P_p)=\lhd _p$.

(d)  If $\langle q_1,q_2 \rangle\in \lhd$ and $p\in \mathcal G$
where $\langle q_1,q_2 \rangle\in\lhd_p$, then by (iii),
$q_1<_{\mathbb Q}q_2$. \hfill $\Box$
\begin{cla}\rm\label{T5035}
(a) $\langle A ,\vartriangleleft \rangle$ is a random poset, for
each $x\in (-\infty , \infty ]$ and each set $A$  satisfying
\begin{equation}\label{EQ5012}
(-\infty ,x)\cap J \subset  A\subset (-\infty , x)\cap \mathbb Q
\end{equation}

(b) If $J \subset A \subset \mathbb Q$ then $\langle A ,\lhd
\rangle$ (in particular, $\langle \mathbb Q ,\vartriangleleft
\rangle$) is a random poset.

(c)  If $C \subset {\mathbb Q}$ and $\max _{\langle \mathbb Q ,
<_\mathbb Q \rangle}C$ exists, then $\langle C  ,\lhd \rangle$ is
not a random poset.
\end{cla}
\dok (a) By Claim \ref{T5034}(b), $\langle A ,\vartriangleleft
\rangle$ is a strict partial order. Let $\langle L,G,U \rangle \in
C(A ,\lhd)$. Then
\begin{equation}\label{EQ5013}
L\cup G \cup U \subset A \;\;\land \;\; L\cap G=G\cap U=U\cap A =\emptyset,
\end{equation}
\begin{equation}\label{EQ5014}
\forall l\in L \;\; \forall g\in G \;\; \forall u\in U \;\;
(\langle l,g \rangle \in \lhd \;\land\; \langle u,l \rangle
\not\in \lhd \;\land\; \langle g,u \rangle \not\in \lhd ).
\end{equation}
We show that $\langle A, \lhd \rangle _{\langle L,G,U \rangle}\neq
\emptyset$. For $\langle L,G,U \rangle\neq \langle
\emptyset,\emptyset,\emptyset \rangle$ we have two cases.

{\it Case 1:} $G\neq \emptyset$. Let $p\in \mathcal G \cap
\mathcal D _{\langle L,G,U \rangle ,1}$. Then
\begin{equation}\label{EQ5015}
L\cup G \cup U \subset P_p .
\end{equation}
First we show that $\langle L,G,U \rangle \in C(p)$. Let $l\in L$,
$g\in G$ and $u\in U$. By (\ref{EQ5014}), (\ref{EQ5015}) and Claim
\ref{T5034}(c) we have $\langle l,g \rangle\in \lhd _p$ and (C1)
is true. Since $\lhd _p \subset \lhd$ by (\ref{EQ5014}) we have
$\langle u,l \rangle \not\in \lhd _p$ and $\langle g,u \rangle
\not\in \lhd _p$ and (C2) and (C3) are true as well.

Since $p\in \mathcal D _{\langle L,G,U \rangle ,1}$ there is $q\in
p_{\langle L,G,U \rangle}\cap J$. We prove that $q\in \langle A,
\lhd \rangle _{\langle L,G,U \rangle}$.
 For a $g\in G$ we have $q\lhd _p g$
and, by (iii), $q <_\mathbb Q g$. By (\ref{EQ5012}) and
(\ref{EQ5013}) we have $g\in G \subset A \subset (-\infty , x)$
and, hence $q<_\mathbb Q g <_\mathbb R x$, thus $q\in (-\infty ,
x)\cap J\subset A$. Let $l\in L$, $g\in G$ and $u\in U$. Since
$q\in p_{\langle L,G,U \rangle}$ we have $l \lhd _p q \lhd _p g$
and $\lhd _p\subset \lhd $ implies $l \lhd  q \lhd  g$. Thus (S1)
and (S2) are true. Suppose that $\neg\; q \parallel _{\langle A,
\lhd \rangle} u$. Since $q\not\in U$ we have $q\neq u$ and, hence,
$q\lhd u$ or $u\lhd q$. But then, since $u,q\in P_p$, by  Claim
\ref{T5034}(c) we would have $q\lhd _p u$ or $u\lhd _p q$, which
is impossible because $q\in p_{\langle L,G,U \rangle}$. So (S3) is
true as well.

{\it Case 2:} $G= \emptyset$. By (\ref{EQ5012}) and (\ref{EQ5013})
we have $L\cup G \cup U\subset (-\infty ,x)$, which implies
$m_{\langle L,G,U \rangle}<x$ and, hence, there is $m\in \mathbb
N$ such that
\begin{equation}\label{EQ5017} \textstyle
m_{\langle L,G,U \rangle} +\frac{1}{m}<x .
\end{equation}
Let $p\in \mathcal G \cap \mathcal D _{\langle L,G,U \rangle ,m}$.
Then (\ref{EQ5015}) holds again and exactly like in Case 1 we show
that $\langle L,G,U \rangle \in C(p)$. Thus, since  $p\in \mathcal
D _{\langle L,G,U \rangle ,m}$ there is $q\in p_{\langle L,G,U
\rangle}\cap (m_{\langle L,G,U \rangle} , m_{\langle L,G,U
\rangle} +\frac{1}{m})\cap J$ and, by (\ref{EQ5017}), $q\in J\cap
(-\infty , x)$. Thus, by (\ref{EQ5012}), $q\in A$ and exactly like
in Case 1 we prove that $q\in \langle A, \lhd \rangle _{\langle
L,G,U \rangle}$.

(b) Follows from (a) for $x=\infty$.

(c) Suppose that $\max _{\langle \mathbb Q , <_\mathbb Q
\rangle}C=q$ and that $\langle C  ,\lhd \rangle$ is a random
poset. Then $C_{\langle \{q\},\emptyset,\emptyset
\rangle}\neq\emptyset$ and, by (S1), there is $q_1\in C$ such that
$q\lhd q_1$, which, by Claim \ref{T5034}(d) implies  $q<_\mathbb Q
q_1$. A contradiction with the maximality of $q$.
\kdok

\noindent
For $y \in M$ let us take $I_y \in [J_y \cap (-\infty,y)]^{|L_y|-1}$
and define $A_{-\infty }= \emptyset$ and
$$ \textstyle
A_x=\big(J \cap (-\infty,x)  \big)\cup \bigcup{_{y \in M \cap
(-\infty,x)}}I_y, \quad \mbox{for }  x \in (-\infty,\infty];
$$
\vspace{-6mm}
$$\textstyle
A_x^+=A_x\cup I_x, \quad \mbox{for } x \in M.
$$
Since $J \subset A_\infty ^+  \subset \mathbb Q$, by Claim
\ref{T5035}(b) $\langle A_\infty ^+ , \lhd \rangle $ is a random
poset and we construct a maximal chain ${\mathcal L}$ in $\langle
\mathbb P(A_\infty ^+ , \lhd ), \subset\rangle$, such that
${\mathcal L}\cong L$.
\begin{cla}  \rm \label{T5026}
The sets  $ A_x $, $x\in [-\infty,\infty]$ and $A_x^{+}$, $x \in M$ are subsets of the set $A_{\infty}^+$ and of $\mathbb Q$.
In addition, for each $x, x_1, x_2 \in [-\infty,\infty]$ we have

(a) $A_x \subset (-\infty,x)$;

(b) $A_x^+ \subset (-\infty,x)$, if $x \in M$;

(c) $x_1 < x_2 \Rightarrow A_{x_1}\varsubsetneq A_{x_2}$;

(d) $M \ni x_1 < x_2 \Rightarrow  A_{x_1}^+ \varsubsetneq A_{x_2}$;

(e)  $|A_x^+\setminus A_x|=|L_x|-1$, if $x \in M$;

(f) $ A_x \in \mathbb P (A_{\infty}^+ )$, for each $x\in (-\infty
, \infty ]$.

(g) $ A_x ^+ \in \mathbb P (A_{\infty}^+  )$ and $[A_x , A_x
^+]_{\mathbb P (A_{\infty}^+  )} = [A_x , A_x ^+]_{P( A_x ^+ )}$,
for each $x\in M$.
\end{cla}
\dok Statements (c) and (d) are true since $J$ is a dense subset
of $\mathbb Q$; (a), (b) and (e) follow from the definitions of
$A_x$ and $A_x^+$ and the choice of the sets $I_y$. Since $J \cap
(-\infty,x) \subset A_x \subset A_x^+ \subset (-\infty,x)\cap
\mathbb Q$, (f) and (g) follow from Claim \ref{T5035}(a).
\kdok

\noindent Now, for $x \in [-\infty,\infty]$ we define chains
${\mathcal L}_x \subset \mathbb P (A_{\infty}^+  )\cup \{
\emptyset \}$ in the following way.

For $x \not\in M$ we define ${\mathcal L}_x=\{A_x\}$. In particular, ${\mathcal L}_{-\infty }=\{ \emptyset \}$.

For $x \in M$, using Claim \ref{T5026} and Lemma \ref{T5028} we obtain a set ${\mathcal L}_x \subset [A_x , A_x ^+]_{P(A_x ^+ )}$
such that $\langle \mathcal L _x , \varsubsetneq \rangle \cong \langle L_x , <_x \rangle $ and
\begin{equation}\label{EQ5002}\textstyle
A_x, A_x^+ \in {\mathcal L}_x \subset [A_x,A_x^+]_{\mathbb P
(A_{\infty}^+ )},
\end{equation}
\begin{equation}\label{EQ5003}\textstyle
\bigcup {\mathcal  A},\bigcap {\mathcal  B}\in {\mathcal  L}_x \;\;\mbox{ and }\;\;
|\bigcap {\mathcal  B}\setminus \bigcup {\mathcal  A}|\leq 1,
\mbox{ for each cut }\langle {\mathcal A},{\mathcal B}\rangle \mbox{ in }{\mathcal L}_x.
\end{equation}
For ${\mathcal A},{\mathcal B} \subset \mathbb P (A_{\infty}^+)$
we will write ${\mathcal A}\prec {\mathcal B}$ iff $A\varsubsetneq
B$, for each $A \in {\mathcal A}$ and $B \in {\mathcal B}$.
\begin{cla}  \rm \label{T5029}
Let ${\mathcal L}=\bigcup_{x \in [-\infty,\infty]}{\mathcal L}_x$. Then

(a) If $-\infty \leq x_1 <x_2  \leq \infty $, then ${\mathcal L}_{x_1}\prec {\mathcal L}_{x_2}$
and $\bigcup{\mathcal L}_{x_1} \subset A_{x_2}\subset \bigcup{\mathcal L}_{x_2}.$

(b) ${\mathcal L}$ is a chain in $\langle \mathbb P
(A_{\infty}^+)\cup \{ \emptyset \},\subset\rangle$ isomorphic to
$L= \sum_{x \in [-\infty,\infty]}L_x$.

(c) ${\mathcal L}$ is a maximal chain in $\langle \mathbb P
(A_{\infty}^+ )\cup \{ \emptyset \},\subset\rangle$.
\end{cla}
\dok
(a) Let $A \in {\mathcal L}_{x_1}$ and $B \in {\mathcal L}_{x_2}$.
If $x_1 \in (-\infty,\infty]\setminus M$, then,
by (\ref{EQ5002}) and Claim \ref{T5026}(c) we have $A=A_{x_1}\varsubsetneq A_{x_2} \subset B$.
If $x_1 \in M$, then, by (\ref{EQ5002}) and Claim \ref{T5026}(d),
$A \subset A_{x_1}^+\varsubsetneq A_{x_2} \subset B$. The second statement
follows from $A_{x_2} \in {\mathcal L}_{x_2}$.

(b) By (a), $\langle [-\infty,\infty],<\rangle \cong \langle \{{\mathcal L}_x:x \in [-\infty,\infty]\},\prec \rangle$.
Since ${\mathcal L}_x \cong L_x$, for $x \in [-\infty,\infty]$, we have
$\langle {\mathcal L},\varsubsetneq\rangle
\cong\sum_{x \in [-\infty,\infty]}\langle {\mathcal L}_x,\varsubsetneq\rangle
\cong \sum_{x \in [-\infty,\infty]}L_x =L$.

(c) Suppose that $C\in \mathbb P (A_{\infty}^+ )\cup \{ \emptyset
\}$ witnesses that ${\mathcal L}$ is not maximal. Clearly
${\mathcal L}={\mathcal A}\dot{\cup } {\mathcal B}$ and ${\mathcal
A}\prec {\mathcal B}$, where ${\mathcal A}=\{ A \in {\mathcal L}:
A \varsubsetneq C \}$ and  ${\mathcal B}=\{ B \in {\mathcal L}: C
\varsubsetneq B \}$. Now $\emptyset \in {\mathcal L}_{-\infty }$
and, since $\infty \in M$, by (\ref{EQ5002}) we have $A_{\infty}^+
\in {\mathcal L}_{\infty}$. Thus $\emptyset ,A_{\infty}^+ \in
{\mathcal L}$, which implies ${\mathcal A},{\mathcal B}\neq
\emptyset$ and, hence, $\langle {\mathcal A},{\mathcal B}\rangle$
is a cut in $\langle {\mathcal L},\varsubsetneq\rangle$. By
(\ref{EQ5002}) we have  $\{ A_x :x \in (-\infty,\infty ]\}\subset
{\mathcal L}\setminus \{ \emptyset \}$ and, by Claim
\ref{T5026}(a), $\bigcap ({\mathcal L}\setminus \{ \emptyset
\})\subset \bigcap_{x \in (-\infty,\infty]}A_x \subset \bigcap
_{x\in (-\infty,\infty]}(-\infty,x)=\emptyset$, which implies
${\mathcal A}\neq \{ \emptyset\}$. Clearly,
\begin{equation}\label{EQ5004}\textstyle
\bigcup{\mathcal A} \subset C \subset \bigcap {\mathcal B}.
\end{equation}
\noindent {\it Case 1:}
${\mathcal A}\cap {\mathcal L}_{x_0}\neq\emptyset$ and ${\mathcal B}\cap {\mathcal L}_{x_0}\neq \emptyset$,
for some $x_0 \in (-\infty,\infty]$.
Then $|{\mathcal L}_{x_0}|>1$, $x_0 \in M$ and
$\langle {\mathcal A}\cap{\mathcal L}_{x_0},{\mathcal B}\cap {\mathcal L}_{x_0}\rangle$
is a cut in ${\mathcal L}_{x_0}$ satisfying (\ref{EQ5003}). By (a),
${\mathcal A}=\bigcup_{x<x_0}{\mathcal L}_x\cup ({\mathcal A}\cap
{\mathcal L}_{x_0})$ and, consequently, $\bigcup{\mathcal A}=\bigcup({\mathcal A}\cap{\mathcal L}_{x_0}) \in {\mathcal L}$. Similarly,
$\bigcap {\mathcal B}=\bigcap({\mathcal B}\cap{\mathcal L}_{x_0})\in {\mathcal L}$ and, since $|\bigcap {\mathcal B}\setminus \bigcup{\mathcal A}|\leq 1$,
by (\ref{EQ5004}) we have $C\in {\mathcal L}$. A contradiction.

\vspace{2mm}
\noindent
{\it Case 2:} $\neg$ Case 1.
Then for each $x \in (-\infty,\infty]$ we have ${\mathcal L}_x \subset {\mathcal A}$
or ${\mathcal L}_x \subset {\mathcal B}$. Since
${\mathcal L}={\mathcal A}\;{\buildrel {.}\over{\cup}}\; {\mathcal B}$, ${\mathcal A}\neq \{ \emptyset\}$
and ${\mathcal A}, {\mathcal B}\neq\emptyset$, the sets
${\mathcal A}'=\{x \in (-\infty,\infty]:{\mathcal L}_x\subset {\mathcal A}\}$ and
${\mathcal B}'=\{x \in(-\infty,\infty]:{\mathcal L}_x \subset {\mathcal B}\}$
are non-empty and $(-\infty,\infty]={\mathcal A}'\;{\buildrel {.} \over {\cup}}\;{\mathcal B}'$.
Since ${\mathcal A}\prec {\mathcal B}$, for $x_1 \in {\mathcal A}'$ and $x_2 \in {\mathcal B}'$ we have
${\mathcal L}_{x_1}\prec {\mathcal L}_{x_2}$ so, by (a), $x_1<x_2$.
Thus $\langle {\mathcal A}',{\mathcal B}' \rangle$ is a cut in $(-\infty,\infty]$ and,
consequently, there is $x_0 \in (-\infty,\infty]$ such that
$x_0=\max {\mathcal A}'$ or $x_0=\min {\mathcal B}'$.

\vspace{2mm}
\noindent
{\it Subcase 2.1:} $x_0=\max {\mathcal A}'$.
Then $x_0<\infty$ because $\mathcal B \neq \emptyset$ and
${\mathcal A}=\bigcup_{x\leq x_0}{\mathcal L}_x$ so, by (a),
$\bigcup{\mathcal A}=\bigcup_{x\leq x_0}\bigcup{\mathcal
L}_x=\bigcup_{x<x_0}\bigcup{\mathcal L}_x \cup \bigcup{\mathcal
L}_{x_0}=\bigcup{\mathcal L}_{x_0}$ which, together with
(\ref{EQ5002}) implies
%$$ \textstyle
\begin{equation}\label{EQ5005}\textstyle
\bigcup{\mathcal A}=\left\{\begin{array}{ll}
                         A_{x_0} & \mbox{if } x_0 \not\in M, \\
                         A_{x_0}^+ & \mbox{if } x_0 \in M.
                         \end{array}
                         \right.
\end{equation}
%$$
Since ${\mathcal B}=\bigcup_{x\in (x_0,\infty]}{\mathcal L}_x$, we
have $\bigcap {\mathcal B}=\bigcap_{x \in (x_0,\infty]}\bigcap{\mathcal L}_x$. By (\ref{EQ5002})
$\bigcap{\mathcal L}_x=A_x$, so we have
$
\bigcap {\mathcal B}  =  \big(\bigcap{_{x \in (x_0,\infty]}}(-\infty,x)\cap J\big)
                                      \cup \big(\bigcap{_{x \in (x_0,\infty]}}\bigcup{_{y \in M\cap(-\infty,x)}}I_y\big)
                      =  \big((-\infty,x_0]\cap J \big) \cup \bigcup{_{y \in M \cap(-\infty,x_0]}}I_y
                      =  \textstyle A_{x_0} \cup \big(\{x_0\}\cap J \big)
                                      \cup \bigcup{_{y \in M\cap\{x_0\}}}I_y ,
$
so
%$$ \textstyle
\begin{equation}\label{EQ5006}\textstyle
\bigcap{\mathcal B}=\left\{\begin{array}{llr}
                          A_{x_0} & \mbox{if} & x_0 \notin J \quad \wedge \quad x_0 \notin M, \\
                          A_{x_0}\cup\{x_0\} & \mbox{if} & x_0 \in J \quad \wedge \quad x_0 \notin M, \\
                          A_{x_0}^+ & \mbox{if} & x_0 \notin J \quad \wedge \quad x_0 \in M, \\
                          A_{x_0}^+ \cup \{x_0\} & \mbox{if} & x_0 \in J \quad \wedge \quad x_0 \in M.
                          \end{array}
                          \right.
\end{equation}
%$$
If $x_0 \not\in J$, then, by (\ref{EQ5004}), (\ref{EQ5005}) and (\ref{EQ5006}), we have $ \bigcup \mathcal A =\bigcap \mathcal B =C \in \mathcal L$.
A contradiction.

If $x_0 \in J$ and $x_0 \not\in M$, then $\bigcup \mathcal A =
A_{x_0}$ and $\bigcap \mathcal B = A_{x_0} \cup \{ x_0\}$. So, by
(\ref{EQ5004}) and since $C\not\in \mathcal L$ we have $C=\bigcap
\mathcal B$. But, by Claim \ref{T5026}(a), $x_0 = \max \bigcap
\mathcal B$ so, by Claim \ref{T5035}(c), $C \not\in \mathbb P
(A_\infty ^+ )$. A contradiction.

If $x_0 \in J$ and $x_0 \in M$, then $\bigcup \mathcal A =
A_{x_0}^+$ and $\bigcap \mathcal B = A_{x_0}^+ \cup \{ x_0\}$.
Again, by (\ref{EQ5004}) and since $C\not\in \mathcal L$ we have
$C=\bigcap \mathcal B$. By Claim \ref{T5026}(b), $x_0 = \max
\bigcap \mathcal B$ so, by Claim \ref{T5035}(c), $C \not\in \mathbb P (A_\infty ^+ )$. A contradiction.

\vspace{2mm}
\noindent
{\it Subcase 2.2:} $x_0=\min {\mathcal B}'$.
Then, by (\ref{EQ5002}), $A_{x_0} \in {\mathcal L}_{x_0} \subset {\mathcal B}$
which, by (a), implies $\bigcap{\mathcal B}=A_{x_0}$. Since $A_x \in {\mathcal L}_x$, for $x \in (-\infty,\infty]$ and
${\mathcal A}=\bigcup_{x<x_0}{\mathcal L}_x$ we have
$\bigcup{\mathcal A}
=\bigcup_{x<x_0}\bigcup{\mathcal L}_x \supset \bigcup_{x<x_0}A_x
=\bigcup_{x<x_0}\big((-\infty,x)\cap J \big)\cup \bigcup_{x<x_0}\bigcup_{y \in M \cap(-\infty,x)}I_y$
$=\big((-\infty,x_0)\cap J \big)\cup \bigcup_{y \in M \cap(-\infty,x_0)}I_y
=A_{x_0}$
so $A_{x_0} \subset \bigcup{\mathcal A}\subset \bigcap{\mathcal B} = A_{x_0}$,
which implies $C=A_{x_0} \in {\mathcal L}$. A contradiction.
\kdok

\noindent {\bf Case II:} $- \infty \not\in M \not\ni \infty$. Then
$L_\infty =\{ \max L\}$ and the sum $L+1$ belongs to Case I. So,
there are a maximal chain ${\mathcal L}$ in $\langle \mathbb P
(\mathbb D )\cup \{\emptyset \}, \subset \rangle $ and an
isomorphism $f: \langle L+1 , < \rangle \rightarrow \langle
{\mathcal L} , \subset \rangle $. Then $A=f(\max L)\in \mathbb P
(\mathbb D )$ and ${\mathcal L}'=f[L]\cong L$. By the maximality
of ${\mathcal L}$, ${\mathcal L}'$ is a maximal chain in $\langle
\mathbb P (A)\cup \{\emptyset \}, \subset \rangle $.

\vspace{2mm}

\noindent
{\bf Case III:} $- \infty \in M$. Then $L=\sum_{x \in [-\infty,\infty]}L_x$, (L1) and (L2) of Claim \ref{T5024} hold and

(L3$'$) $L_{-\infty}$ is a countable complete linear order with $0_{L_{-\infty}}$ non-isolated.

\noindent Clearly $L=L_{-\infty }+ L^+$, where $L^+ =\sum_{x \in
(-\infty,\infty]}L_x = \sum_{y \in (0,\infty]}L_{\ln y}$ (here
$\ln \infty =\infty$). Let $L_y'$, $y \in [-\infty,\infty]$, be
disjoint linear orders such that $L_y'\cong 1$, for $y \in
[-\infty , 0]$, and $L_y'\cong L_{\ln y}$, for $y\in (0, \infty
]$. Now $\sum_{y \in [-\infty,\infty]}L_y' \cong [-\infty , 0] +
L^+$ belongs to Case I or Case II and we obtain a maximal chain
${\mathcal L}$ in $\mathbb P (\mathbb D )\cup \{ \emptyset \}$ and
an isomorphism $f: \langle [-\infty , 0] + L^+ , < \rangle
\rightarrow \langle {\mathcal L}, \subset \rangle$. Clearly, for
$A_0=f(0)$ and ${\mathcal L}^+=f[L^+]$ we have $A_0\in {\mathcal
L}$ and ${\mathcal L}^+\cong L^+$.

By (L3$'$) and the fact that (b) $\Rightarrow$ (a) for countable
$L$'s, $\mathbb P (A_0)\cup \{ \emptyset \}$ contains a maximal
chain ${\mathcal L}_{-\infty}\cong L_{-\infty }$. Clearly $A_0\in
{\mathcal L}_{-\infty}$ and ${\mathcal L}_{-\infty} \cup {\mathcal
L}^+ \cong L_{-\infty }+L^+ =L$. Suppose that $B$ witnesses that
${\mathcal L}_{-\infty} \cup {\mathcal L}^+ $ is not a maximal
chain in $\mathbb P (\mathbb D )\cup \{ \emptyset \}$. Then either
$A_0 \varsubsetneq B$, which is impossible since ${\mathcal L}$ is
maximal in $\mathbb P (\mathbb D )\cup \{ \emptyset \}$, or
$B\varsubsetneq A_0$, which is impossible since ${\mathcal
L}_{-\infty}$ is maximal in $\mathbb P(A_0)\cup \{ \emptyset \}$.
\hfill $\Box$
\section{Maximal chains in $\mathbb P(\mathbb B_n)$}\label{S6}
\begin{te}  \rm \label{T5994}
For $n\in\mathbb N$ and each $\mathbb R$-embeddable complete
linear order $L$ with $0_L$ non-isolated there is a maximal chain
in $\langle \mathbb P(\mathbb B_n)\cup \{\emptyset\},\subset
\rangle$ isomorphic to $L$.
\end{te}
\dok Let the order on $\mathbb B_n=\bigcup_{i<n}\mathbb
Q_i=\bigcup_{i<n}\{i\}\times \mathbb Q$ be given by
$$\textstyle
\langle i_1,q_1 \rangle<\langle i_2,q_2 \rangle\Leftrightarrow
i_1=i_2\ \wedge\ q_1<_\mathbb Q q_2.
$$
Clearly, $\langle \mathbb Q,<_{\mathbb Q} \rangle\cong_{f_i}
\langle \mathbb Q_i,< \rangle$, where $f_i(q)=\langle i,q
\rangle$, for all $q\in \mathbb Q$ and, hence, $\mathbb P(\mathbb
Q_i)=\{\{i\}\times C: C\in \mathbb P(\mathbb Q)\}$. If $f: \mathbb
B_n \hookrightarrow \mathbb B_n$, then for each $i<n$ the
restriction $f| \mathbb Q_i$ is an isomorphism, thus there is $j_i
< n$ such that $f[\mathbb Q_i ]\subset \mathbb Q _{j_i}$ and,
moreover, $f[\mathbb Q_i ]\in \mathbb P (\mathbb Q _{j_i})$.
Clearly, $i_1\neq i_2$ implies $j_{i_1}\neq j_{i_2}$ and, thus, we
have
\begin{equation}\label{EQ5994}\textstyle
\mathbb P(\mathbb B_n)=\big\{\bigcup_{i<n}\{i\}\times C_i: \forall
i<n\;\; C_i\in \mathbb P(\mathbb Q)\big\}.
\end{equation}
Now, by Theorem 6 of \cite{K1}, there is a maximal chain $\mathcal
L$ in $\langle \mathbb P(\mathbb Q)\cup \{\emptyset\},\subset
\rangle$ isomorphic to $L$. For $A\in \mathcal L\setminus
\{\emptyset\}$ let
\begin{equation}\label{EQ5993}\textstyle
A^*=( \{0\}\times A ) \cup \bigcup_{0<i<n}\{i\}\times \mathbb Q.
\end{equation}
By (\ref{EQ5994}) we have $\mathcal L^*=\{A^*:A\in \mathcal
L\setminus\{\emptyset\}\}\cup \{\emptyset\}\subset \mathbb
P(\mathbb B_n)\cup \{\emptyset\}$ and, clearly, $\langle \mathcal
L^*,\subset \rangle$ is a chain in $\langle \mathbb P(\mathbb
B_n)\cup \{\emptyset\},\subset \rangle$ isomorphic to $\langle
\mathcal L,\subset \rangle$ and, hence, to $L$. Suppose that some
$C=\bigcup_{i<n}\{i\}\times C_i\in \mathbb P(\mathbb B_n)$
witnesses that $\mathcal L^*$ is not a maximal chain. By
(\ref{EQ5994}) and (\ref{EQ5993}) $C\subset \bigcap _{A\in
\mathcal L\setminus \{\emptyset\}}A^*$ would imply $\mathbb P
(\mathbb Q )\ni C_0\subset \bigcap (\mathcal L\setminus
\{\emptyset\})$, which is impossible ($\mathcal L$ is a maximal
chain in $\mathbb P(\mathbb Q)\cup \{ \emptyset \}$ and
$C_0\setminus F\in \mathbb P (\mathbb Q )$ for each finite
$F\subset C_0$). Thus there is $A\in \mathcal L\setminus
\{\emptyset\}$ such that $A^*\subset C$ and, by (\ref{EQ5993}),
\begin{equation}\label{EQ5992}\textstyle
C=\{0\}\times C_0\cup \bigcup_{0<i<n}\{i\}\times \mathbb Q.
\end{equation}
Since $\mathcal L^*\cup \{C\}$ is a chain, for each $A\in \mathcal
L\setminus \{\emptyset\}$ we have $A^*\subsetneq C\ \vee\
C\subsetneq A^*$ which together with (\ref{EQ5993}) and
(\ref{EQ5992}) implies $A\subsetneq C_0$ or $C_0\subsetneq A$. A
contradiction to the maximality of $\mathcal L$. \hfill $\Box$
\begin{te}  \rm \label{T5993}
For each $\mathbb R$-embeddable complete linear order $L$ with
$0_L$ non-isolated there is a maximal chain in $\langle \mathbb
P(\mathbb B_{\omega})\cup\{\emptyset\},\subset \rangle$ isomorphic
to $L$.
\end{te}
\dok Let $x_0=\infty$, let $\langle x_n:n\in\mathbb N \rangle$ be a
descending sequence in $\mathbb R\setminus \mathbb Q$ without a
lower bound and let $\mathbb B_{\omega}=\langle \mathbb Q,<_{\omega}
\rangle=\bigcup_{i\in\omega}\langle (x_{i+1},x_i)\cap \mathbb
Q,<_i \rangle$ where
$$\textstyle
q_1<_{\omega}q_2\Leftrightarrow \exists i\in \omega \;\;
(q_1,q_2\in (x_{i+1},x_i)\ \wedge\ q_1<_\mathbb Q q_2 ).
$$
Then for the sets $\mathbb Q_i=(x_{i+1},x_i)\cap \mathbb Q$, $i\in
\omega$, we have $\langle \mathbb Q_i,<_i \rangle \cong \langle
\mathbb Q,<_{\mathbb Q} \rangle$, which implies $\mathbb P(\mathbb
Q_i,<_i)\cong \mathbb P(\mathbb Q,<_{\mathbb Q})$. As in the proof
of Theorem \ref{T5994} we obtain
\begin{equation}\label{EQ5991}\textstyle
\mathbb P(\mathbb B_{\omega})=\{\bigcup_{i\in S}C_i:S\in
[\omega]^{\omega}\ \wedge\ \forall i\in S\;\; C_i\in \mathbb
P(\mathbb Q_i) \}.
\end{equation}
Let $L$ be a linear order with the given properties and, first,
let $|L|=\omega$. Clearly the family $\dense(\mathbb Q_i)$ of
dense subsets of $\mathbb Q _i$ is a subset of $\mathbb P(\mathbb
Q_i)$ and by (\ref{EQ5991}) we have $\textstyle \mathcal
P=\big\{\bigcup_{i\in \omega}C_i: \forall i\in \omega \;\; C_i\in
\dense(\mathbb Q_i)\big\}\subset\mathbb P(\mathbb B_{\omega}). $
It is easy to check that $\mathcal P$ is a positive family on
$\mathbb Q$ so, by Theorem \ref{T5037}(b), there is a maximal
chain in $\langle \mathbb P(\mathbb
B_{\omega})\cup\{\emptyset\},\subset \rangle$ isomorphic to $L$.

Now, let $|L|>\omega$. Then, by Claim \ref{T5024}, we can assume
that $L=\sum_{x\in [-\infty,\infty]}L_x$, where conditions (L1-L3)
from Claim \ref{T5024} are satisfied. We distinguish two cases.

{\it Case 1:} $-\infty\notin M$. Then, by the construction from
\cite{K1} (if $(0,1]$ is replaced by $(-\infty,\infty]$ and
$A_1^+$ by $\mathbb Q$), there is a maximal chain $\mathcal L$ in
$\langle \mathbb P(\mathbb Q)\cup \{\emptyset\}, \subset \rangle$
such that
\begin{equation}\label{EQ5990}
\forall A\in \mathcal L\setminus \{\emptyset\}\ \exists x\in
(-\infty,\infty]\ \big(A\subset (-\infty,x)\ \wedge\ A\ \mbox{is
dense in}\ (-\infty,x)\big)
\end{equation}
and $\mathcal L \cong L$. Now we prove
\begin{equation}\label{EQ5989}
\mathcal L\setminus\{\emptyset\}\subset \mathbb P(\mathbb
B_{\omega})\subset \mathbb P(\mathbb Q,<_{\mathbb Q}).
\end{equation}
Let $A\in \mathcal L\setminus \{\emptyset\}$, let $x$ be the real
corresponding to $A$ in the sense of (\ref{EQ5990}) and let
$i_0=\min\{i\in \omega: (-\infty,x)\cap
(x_{i+1},x_i)\neq\emptyset\}$. Then $x_{i_0+1}<x\leq x_{i_0}$ and,
by (\ref{EQ5990}) the set $C_{i_0}=A\cap (x_{i_0+1},x)$ is dense
in $(x_{i_0+1},x)$ and, hence, $C_{i_0}\in \mathbb P(\mathbb
Q_{i_0})$. Similarly,  $C_i=A\cap (x_{i+1},x_i)\in \mathbb
P(\mathbb Q_i)$, for all $i>i_0$. Since $A\subset \mathbb Q$, we
have $A=\bigcup_{i\geq i_0}C_i$ and, by (\ref{EQ5991}), $A\in
\mathbb P(\mathbb B_{\omega})$. So the first inclusion of
(\ref{EQ5989}) is proved.

Let $C=\bigcup_{i\in S}C_i\in \mathbb P(\mathbb B_{\omega})$. By
(\ref{EQ5991}) for each $i\in S$ we have $C_i\cong \mathbb Q _i
\cong \mathbb Q$ and, hence, $C\cong \sum_{\omega^*}\mathbb
Q\cong\mathbb Q$. The second inclusion of (\ref{EQ5989}) is proved
as well.

By (\ref{EQ5989}) we have $\mathcal L\subset \mathbb P(\mathbb
B_{\omega})\cup \{\emptyset\}\subset \mathbb P(\mathbb
Q,<_{\mathbb Q})\cup \{\emptyset\}$ and, clearly, $\mathcal L$ is
a chain in $\mathbb P(\mathbb B_{\omega})\cup \{\emptyset\}$.
Suppose that $\mathcal L\cup \{C\}$ is a chain, for some $C\in
(\mathbb P(\mathbb B_{\omega})\cup\{\emptyset\})\setminus \mathcal
L$. Then, by (\ref{EQ5989}), $C\in \mathbb P(\mathbb Q,<_{\mathbb
Q})$ and $\mathcal L$ would not be a maximal chain in the poset
$\langle \mathbb P(\mathbb Q,<_{\mathbb
Q})\cup\{\emptyset\},\subset \rangle$. So $\mathcal L $ is a
maximal chain in $\langle \mathbb P(\mathbb
B_{\omega})\cup\{\emptyset\},\subset \rangle$ and $\mathcal L\cong
L$.

{\it Case 2:} $-\infty\in M$. Then we proceed as in (III) of the proof of Theorem \ref{T5041}.
\hfill $\Box$
\section{Maximal chains in $\mathbb P(\mathbb C_n)$}\label{S7}
\begin{te} \rm \label{T5992}
For all $n\in \mathbb N$ and each $\mathbb R$-embeddable complete
linear order $L$ with $0_L$ non-isolated there is a maximal chain
in $\langle \mathbb P(\mathbb C_n)\cup\{\emptyset\}, \subset
\rangle$ isomorphic to $L$.
\end{te}
\dok Let the order $<$ on $\mathbb C _n=\mathbb Q\times n$ be
given by $\langle q_1,i_1 \rangle<\langle q_2,i_2
\rangle\Leftrightarrow q_1<_{\mathbb Q}q_2$. Clearly, the
incomparability relation $a\| b\Leftrightarrow a\nless b\ \wedge\
b\nless a$ on $\mathbb C_n$ is an equivalence relation with the
equivalence classes $\{q\}\times n$, $q\in \mathbb Q$, of size $n$
and the corresponding quotient, $\mathbb C_n/\|$, is isomorphic to
$\langle \mathbb Q,<_{\mathbb Q} \rangle$. Since each element of
$\mathbb P(\mathbb C_n)$ has such classes we have $\mathbb
P(\mathbb C_n)=\{A\times n:A\in \mathbb P(\mathbb Q,<_{\mathbb
Q})\}$. It is easy to see that the mapping $f:\mathbb P(\mathbb
Q,<_{\mathbb Q})\cup\{\emptyset\}\to \mathbb P(\mathbb
C_n)\cup\{\emptyset\}$, given by $f(A)=A\times n$, is an
isomorphism of partial orders $\langle \mathbb P(\mathbb
Q,<_{\mathbb Q})\cup\{\emptyset\}, \subset \rangle$ and $\langle
\mathbb P(\mathbb C_n)\cup\{\emptyset\},\subset \rangle$. Hence
the statement follows from Theorem 6 of \cite{K1}. \hfill $\Box$
\begin{te}  \rm \label{T5991}
For each $\mathbb R$-embeddable complete linear order $L$ with
$0_L$ non-isolated there is a maximal chain in $\langle \mathbb
P(\mathbb C_{\omega})\cup\{\emptyset\},\subset \rangle$ isomorphic
to $L$.
\end{te}
\dok Let the strict order $<$ on $\mathbb C_{\omega}=\mathbb
Q\times\omega =\bigcup_{q\in \mathbb Q}\{q\}\times
\omega=\bigcup_{q\in \mathbb Q}\omega _q$ be given by $\langle
q_1,i_1 \rangle<\langle q_2,i_2 \rangle\Leftrightarrow
q_1<_{\mathbb Q}q_2$. For a set $X\subset \mathbb C_{\omega}$  let
us define $\supp X=\{q\in \mathbb Q : X\cap \omega _q \neq
\emptyset\}$. Now the incomparability classes $\omega _q$ are
infinite and, again, the corresponding quotient, $\mathbb C_\omega
/\|$, is isomorphic to the rational line $\langle \mathbb
Q,<_{\mathbb Q} \rangle$. Since the same holds for the copies of
$\mathbb C _\omega$ it is easy to check that
\begin{equation}\textstyle \label{EQ5022}
\mathbb P(\mathbb C_{\omega})=\{\bigcup_{q\in A}\{q\}\times C_{q}:
A\in \mathbb P(\mathbb Q,<_{\mathbb Q})\ \wedge\ \forall q\in A\
C_{q}\in [\omega]^{\omega}\}.
\end{equation}
\begin{equation}\textstyle \label{EQ5023}
X\subset \mathbb C_{\omega} \land \mbox{ there is }\max\supp X
\Rightarrow X\notin \mathbb P(\mathbb C_{\omega}).
\end{equation}
By (\ref{EQ5022}), $\mathcal P=\{\bigcup_{q\in \mathbb
Q}\{q\}\times C_{q}: \forall q\in \mathbb Q \ C_{q}\in
[\omega]^{\omega}\} \subset\mathbb P(\mathbb C_{\omega})$ and,
clearly, $\mathcal P$ is a positive family so for a countable $L$
the statement follows from Theorem \ref{T5037}(b).

Now, let $L$ be an uncountable linear order.
Then, by Claim \ref{T5024}, we can assume that
$L=\sum_{x\in [-\infty,\infty]}L_x$, where conditions (L1-L3)
from Claim \ref{T5024} are satisfied.

\vspace{2mm} \noindent {\bf Case I:} $-\infty\not\in M\ni \infty$.
Let $\mathbb Q=\bigcup_{y\in M}J_y$ be a partition of $\mathbb Q$
into $|M|$ disjoint dense sets and, for $y \in M$, let $I_y \in
[J_y \cap (-\infty,y)]^{|L_y|-1}$. Let $(-\infty,x)_\mathbb Q=
(-\infty,x)\cap\mathbb Q$ and $\omega ^+=\omega \setminus \{0\}$.
Let us define $A_{-\infty }= \emptyset$ and, for $\ x\in
(-\infty,\infty]$,
$$ \textstyle
A_x=((-\infty,x)_\mathbb Q \times \omega^+)\cup \bigcup_{y \in M
\cap (-\infty,x)}I_y \times\{0\},
$$
\vspace{-6mm}
$$\textstyle
A_x^+=A_x\cup (I_x\times\{0\}), \quad \mbox{for } x \in M.
$$
By (\ref{EQ5022}), $A^+_\infty \cong \mathbb C _\omega$ and we
will construct a maximal chain $\mathcal L \cong L$ in the poset
$\langle \mathbb P(A^+_\infty)\cup\{\emptyset\},\subset \rangle$.
By (\ref{EQ5022}), for each $x\in (-\infty , \infty ]$  and each
set $A\subset \mathbb C _\omega$ we have
\begin{equation}\label{EQ5024}
(-\infty,x)_\mathbb Q \times \omega ^+ \subset A \subset
(-\infty,x)_\mathbb Q \times \omega \;\;\Rightarrow\;\; A\in
\mathbb P (\mathbb C _\omega).
\end{equation}
\begin{cla}  \rm \label{T5998}
The sets  $ A_x $, $x\in [-\infty,\infty]$ and $A_x^{+}$, $x \in
M$ are subsets of the set $A_{\infty}^+$. In addition, for each
$x, x_1, x_2 \in [-\infty,\infty]$ we have

(a) $A_x \subset (-\infty,x)_{\mathbb Q}\times \omega$;

(b) $A_x^+ \subset (-\infty,x)_{\mathbb Q}\times \omega$, if $x
\in M$;

(c) $x_1 < x_2 \Rightarrow A_{x_1}\varsubsetneq A_{x_2}$;

(d) $M \ni x_1 < x_2 \Rightarrow  A_{x_1}^+ \varsubsetneq A_{x_2}$;

(e)  $|A_x^+\setminus A_x|=|L_x|-1$, if $x \in M$;

(f) $ A_x \in \mathbb P (A_{\infty}^+ )$, for each $x\in (-\infty
,\infty ]$.

(g) $ A_x ^+ \in \mathbb P (A_{\infty}^+  )$ and $[A_x , A_x
^+]_{\mathbb P (A_{\infty}^+  )} = [A_x , A_x ^+]_{P( A_x ^+ )}$,
for each $x\in M$.
\end{cla}
\dok Statements (c) and (d) are true since $\mathbb Q$ is a dense
subset of $\mathbb R$; (a), (b) and (e) follow from the
definitions of $A_x$ and $A_x^+$ and the choice of the sets $I_y$.
Since $(-\infty,x)_\mathbb Q \times \omega ^+ \subset A_x \subset
A_x^+ \subset (-\infty,x)_\mathbb Q \times \omega$, (f) and (g)
follow from (\ref{EQ5024}).
\kdok
\noindent Now, for $x \in
[-\infty,\infty]$ we define chains ${\mathcal L}_x \subset \mathbb
P (A_{\infty}^+ )\cup \{ \emptyset \}$ in the following way.

For $x \not\in M$ we define ${\mathcal L}_x=\{A_x\}$. In particular, ${\mathcal L}_{-\infty }=\{ \emptyset \}$.

For $x \in M$, by Claim \ref{T5998} and Lemma \ref{T5028} there is
a set ${\mathcal L}_x \subset [A_x , A_x ^+]_{P(A_x ^+ )}$ such
that $\langle \mathcal L _x ,\varsubsetneq \rangle \cong \langle
L_x , <_x \rangle $ and
\begin{equation}\label{EQ5999}\textstyle
A_x, A_x^+ \in {\mathcal L}_x \subset [A_x,A_x^+]_{\mathbb P
(A_{\infty}^+ )},
\end{equation}
\begin{equation}\label{EQ5998}\textstyle
\bigcup {\mathcal  A},\bigcap {\mathcal  B}\in {\mathcal  L}_x
\;\;\mbox{ and }\;\; |\bigcap {\mathcal  B}\setminus \bigcup
{\mathcal  A}|\leq 1, \mbox{ for each cut }\langle {\mathcal
A},{\mathcal B}\rangle \mbox{ in }{\mathcal L}_x.
\end{equation}
For ${\mathcal A},{\mathcal B} \subset \mathbb P (A_{\infty}^+)$
we will write ${\mathcal A}\prec {\mathcal B}$ iff $A\varsubsetneq
B$, for each $A \in {\mathcal A}$ and $B \in {\mathcal B}$.
\begin{cla}  \rm \label{T5997}
Let ${\mathcal L}=\bigcup_{x \in [-\infty,\infty]}{\mathcal L}_x$. Then

(a) If $-\infty \leq x_1 <x_2  \leq \infty $, then ${\mathcal
L}_{x_1}\prec {\mathcal L}_{x_2}$ and $\bigcup{\mathcal L}_{x_1}
\subset A_{x_2}\subset \bigcup{\mathcal L}_{x_2}.$

(b) ${\mathcal L}$ is a chain in $\langle \mathbb P
(A_{\infty}^+)\cup \{ \emptyset \},\subset\rangle$ isomorphic to
$L= \sum_{x \in [-\infty,\infty]}L_x$.

(c) ${\mathcal L}$ is a maximal chain in $\langle \mathbb P
(A_{\infty}^+ )\cup \{ \emptyset \},\subset\rangle$.
\end{cla}
\dok
The proof of (a) and (b) is a copy of the proof of (a) and (b) of Claim \ref{T5029},
if we replace (\ref{EQ5002}) and Claim \ref{T5026} by (\ref{EQ5999}) and Claim \ref{T5998}.

(c) Suppose that $C\in \mathbb P (A_{\infty}^+ )\cup \{ \emptyset
\}$ witnesses that ${\mathcal L}$ is not maximal. Using
(\ref{EQ5999}) and Claim \ref{T5998}, as in the proof of Claim
\ref{T5029}(c) for ${\mathcal A}=\{ A \in {\mathcal L}: A
\varsubsetneq C \}$ and  ${\mathcal B}=\{ B \in {\mathcal L}: C
\varsubsetneq B \}$ we show that $\langle {\mathcal A},{\mathcal
B}\rangle$ is a cut in $\langle {\mathcal
L},\varsubsetneq\rangle$, ${\mathcal A}\neq \{ \emptyset\}$ and
\begin{equation}\label{EQ5997}\textstyle
\bigcup{\mathcal A} \subset C \subset \bigcap {\mathcal B}.
\end{equation}
\noindent {\it Case 1:} ${\mathcal A}\cap {\mathcal L}_{x_0}\neq\emptyset$ and ${\mathcal B}\cap {\mathcal L}_{x_0}\neq \emptyset$,
for some $x_0 \in (-\infty,\infty]$. Then we obtain
a contradiction exactly like in Claim \ref{T5029}.

\vspace{2mm} \noindent {\it Case 2:} $\neg$ Case 1.
Then like in Claim \ref{T5029} for ${\mathcal A}'=\{x \in (-\infty,\infty]:{\mathcal L}_x\subset {\mathcal A}\}$ and
${\mathcal B}'=\{x \in(-\infty,\infty]:{\mathcal L}_x \subset {\mathcal B}\}$
we show that $\langle {\mathcal A}',{\mathcal B}' \rangle$ is a cut in
$(-\infty,\infty]$. Thus, there is $x_0 \in (-\infty,\infty]$ such that
$x_0=\max {\mathcal A}'$ or $x_0=\min {\mathcal B}'$.

\vspace{2mm} \noindent {\it Subcase 2.1:}
$x_0=\max {\mathcal A}'$. Then like in Claim \ref{T5029} we prove
\begin{equation}\label{EQ5996}\textstyle
\bigcup{\mathcal A}=\left\{\begin{array}{ll}
                         A_{x_0} & \mbox{if } x_0 \not\in M, \\
                         A_{x_0}^+ & \mbox{if } x_0 \in M.
                         \end{array}
                         \right.
\end{equation}
Since ${\mathcal B}=\bigcup_{x\in (x_0,\infty]}{\mathcal L}_x$, we
have $\bigcap {\mathcal B}=\bigcap_{x \in
(x_0,\infty]}\bigcap{\mathcal L}_x$. By (\ref{EQ5999})
$\bigcap{\mathcal L}_x=A_x$, so $ \bigcap {\mathcal B} =
(\bigcap{_{x \in (x_0,\infty]}}(-\infty,x)_{\mathbb Q}\times
\omega^+ )
                                      \cup (\bigcap{_{x \in (x_0,\infty]}}\bigcup{_{y \in M\cap(-\infty,x)}}I_y\times\{0\})
= ((-\infty,x_0]_{\mathbb Q}\times \omega^+ ) \cup \bigcup{_{y \in
M \cap(-\infty,x_0]}}I_y\times\{0\} = A_{x_0} \cup ((\{x_0\}\cap
\mathbb Q )\times \omega ^+) \cup \bigcup{_{y \in
M\cap\{x_0\}}}I_y\times\{0\} , $ so
\begin{equation}\label{EQ5995}\textstyle
\bigcap{\mathcal B}=\left\{\begin{array}{llr}
                          A_{x_0}                                          & \mbox{if} & x_0 \notin \mathbb Q \quad \wedge \quad x_0 \notin M, \\
                          A_{x_0}   \cup (\{x_0\}\times \omega^+ )             & \mbox{if} & x_0 \in \mathbb Q \quad    \wedge \quad x_0 \notin M, \\
                          A_{x_0}^+                                        & \mbox{if} & x_0 \notin \mathbb Q \quad \wedge \quad x_0 \in M, \\
                          A_{x_0}^+ \cup (\{x_0\}\times \omega^+ )             & \mbox{if} & x_0 \in \mathbb Q \quad    \wedge \quad x_0 \in M.
                          \end{array}
                          \right.
\end{equation}
If $x_0 \not\in \mathbb Q$, then, by (\ref{EQ5997}-\ref{EQ5995}),
we have $ \bigcup \mathcal A =\bigcap \mathcal B =C \in \mathcal
L$. A contradiction.

If $x_0 \in \mathbb Q$ and $x_0 \not\in M$, then $\bigcup \mathcal
A = A_{x_0}$ and $\bigcap \mathcal B = A_{x_0} \cup (\{x_0\}\times
\omega ^+ )$. So, by (\ref{EQ5997}) and since $C\not\in \mathcal
L$ we have $C=A_{x_0}\cup S$, where $\emptyset \neq S \subset \{
x_0\}\times \omega ^+$. By Claim \ref{T5998}(a), $x_0 = \max \supp
C$ so, by (\ref{EQ5023}), $C \not\in \mathbb P (A_\infty ^+ )$. A
contradiction.

If $x_0 \in \mathbb Q$ and $x_0 \in M$, then $\bigcup \mathcal A =
A_{x_0}^+$ and $\bigcap \mathcal B = A_{x_0}^+ \cup (\{x_0\}\times
\omega^+)$. Again, by (\ref{EQ5997}) and since $C\not\in \mathcal
L$ we have $C=A_{x_0}\cup S$, where $\emptyset \neq S \subset \{
x_0\}\times \omega ^+$. By Claim \ref{T5998}(b), $x_0 = \max \supp
C$ so, by (\ref{EQ5023}), $C \not\in \mathbb P (A_\infty ^+ )$. A
contradiction.

\vspace{2mm} \noindent {\it Subcase 2.2:} $x_0=\min {\mathcal B}'$.
Then, by (\ref{EQ5999}), $A_{x_0} \in {\mathcal L}_{x_0} \subset {\mathcal B}$ which, by (a), implies
$\bigcap{\mathcal B}=A_{x_0}$. Since $A_x \in {\mathcal L}_x$, for
all $x \in (-\infty,\infty]$ and ${\mathcal A}=\bigcup_{x<x_0}{\mathcal L}_x$ we have
$
\bigcup{\mathcal A}
%=\bigcup_{x<x_0}\bigcup{\mathcal L}_x
\supset \bigcup_{x<x_0}A_x =\bigcup_{x<x_0}((-\infty,x)_{\mathbb
Q}\times \omega^+ )\cup \bigcup_{x<x_0}\bigcup_{y \in
M\cap(-\infty,x)}I_y\times\{0\} =((-\infty,x_0)_{\mathbb Q}\times
\omega ^+ )\cup \bigcup_{y \in M \cap(-\infty,x_0)}I_y\times\{0\}
=A_{x_0} $ so $A_{x_0} \subset \bigcup{\mathcal A}\subset
\bigcap{\mathcal B} = A_{x_0}$, which implies $C=A_{x_0} \in
{\mathcal L}$. A contradiction.
\kdok

\noindent
{\bf Case II:} $-\infty\not\in M\not\ni \infty$ or $-\infty\in M$. Then we proceed like in  Claim \ref{T5029}.
\hfill $\Box$


\footnotesize
\begin{thebibliography}{22}
\bibitem{Cam}
      P.\ J.\ Cameron, D.\ C.\ Lockett,
      Posets, homomorphisms and homogeneity,
      Discrete Math. 310,3 (2010) 604--613.
\bibitem{Day}
      G.\ W.\ Day,
      Maximal chains in atomic Boolean algebras,
      Fund.\ Math.\ 67 (1970) 293-–296.
\bibitem{Fra}
      R.\ Fra\"{\i}ss\'{e},
      Theory of relations, Revised edition, With an appendix by Norbert Sauer,
      Studies in Logic and the Foundations of Mathematics, 145,
      North-Holland, Amsterdam, 2000.
\bibitem{Hodg}
      W. Hodges, Model theory, Encyclopedia of Mathematics and its Applications, 42, Cambridge University Press, Cambridge, 1993.
\bibitem{Kopp}
      S.\ Koppelberg,
      Maximal chains in interval algebras,
      Algebra Univers.\ 27,1 (1990) 32-–43.
\bibitem{Kura}
      K.\ Kuratowski,
      Sur la notion de l'ordre dans la th\'eorie des ensembles,
      Fund.\ Math.\ 2 (1921) 161--171.
\bibitem{K}
       M.\ S.\ Kurili\'c,
       Maximal chains in positive subfamilies of $P(\omega )$,
       Order, 29,1 (2012) 119--129.
\bibitem{K1}
       M.\ S.\ Kurili\'c,
       Maximal chains of copies of the rational line,
       Order,
       in print, %\\
       DOI 10.1007/s11083-012-9273-1
\bibitem{KurTod}
      M.\ S.\ Kurili\'c, S.\ Todor\v cevi\'c,
      Forcing by non-scattered sets,
      Ann.\ Pure Appl.\ Logic 163 (2012) 1299--1308.
\bibitem{Kscatt}
      M.\ S.\ Kurili\'c,
      Posets of copies of countable scattered linear orders,
      submitted.
%\bibitem{Kmaxemb}
%       M.\ S.\ Kurili\'c,
%       Maximally embeddable components,
%       submitted.
%\bibitem{McMo}
%       R.\ McKenzie, J.D.\ Monk,
%       Chains in Boolean algebras,
%       Ann.\ Math.\ Logic 22,2 (1982) 137-–175.
\bibitem{Monk}
       J.\ D.\ Monk,
       Towers and maximal chains in Boolean algebras,
       Algebra Univers.\ 56,3-4 (2007) 337--347.
\bibitem{Rosen}
       J.\ G.\ Rosenstein,
       Linear orderings,
       Pure and Applied Mathematics, 98,
       Academic Press Inc., New York-London, 1982.
\bibitem{Sch}
       J.\ H.\ Schmerl,
       Countable homogeneous partially ordered sets,
       Algebra Univers.\ 9,3 (1979) 317-–321.
\end{thebibliography}
\end{document}